\input amstex
\documentstyle{amsppt}
\document
\NoBlackBoxes
\nologo
\def\p{\partial}
\def\b{\bullet}
\def\D{\Delta}
\def\d{\delta}
\def\A{\Cal{A}}

\centerline{\bf Three constructions of Frobenius manifolds:}

\medskip

\centerline{\bf a comparative study}

\bigskip

\centerline{ Yu.I.Manin}

\medskip

\centerline{\it Max--Planck--Institute f\"ur Mathematik, Bonn}

\bigskip

\centerline{\bf Abstract}

\medskip

The paper studies three classes of Frobenius manifolds:
Quantum Cohomology (topological sigma--models), unfolding spaces of singularities (K.~Saito's theory, Landau--Ginzburg models), and the recent Barannikov--Kontsevich construction starting with the Dolbeault complex
of a Calabi--Yau manifold and conjecturally producing
the $B$--side of the Mirror Conjecture in arbitrary dimension. Each known construction provides the relevant Frobenius manifold with an extra structure
which can be thought of as a version of``non--linear cohomology''. The comparison of thesestructures sheds some light on the general Mirror
Problem: establishing isomorphisms between Frobenius manifolds of different classes. Another theme is the study of tensor products
of Frobenius manifolds, corresponding respectively to
the K\"unneth formula in Quantum Cohomology,
direct sum of singularities in Saito's  theory,
and presumably, the tensor product of the
differential Gerstenhaber--Batalin--Vilkovisky
algebras. We extend the initial Gepner's
construction of mirrors to the context of Frobenius manifolds
and formulate the relevant mathematical conjecture.

\bigskip

\centerline{\bf Contents}

\medskip

\S 0. Introduction

\smallskip

\S 1. Quantum Cohomology and Frobenius manifolds

\smallskip

\S 2. K.~Saito's frameworks

\smallskip

\S 3. Tensor products and their submanifolds

\smallskip

\S 4. Maurer--Cartan and master equations

\smallskip

\S 5. Gerstenhaber--Batalin--Vilkovyski algebras

\smallskip

\S 6. From dGBV--algebras to Frobenius manifolds

\smallskip

References

\newpage

\centerline{\bf \S 0. Introduction}

\medskip

{\bf 0.1. Frobenius manifolds.} Frobenius manifolds 
were introduced and investigated by B. Dubrovin as the axiomatization
of a part of the rich mathematical structure of the Topological
Field Theory (TFT): cf. [D].

\smallskip

 According to
[D] and [M], {\it a Frobenius manifold} is a quadruple 
$(M,\Cal{T}_M^f, g ,A)$.
Here $M$ is a supermanifold in one of the standard categories
($C^{\infty}$, analytic, algebraic, formal, $\dots$),
$\Cal{T}_M^f$  is the sheaf of flat vector fields
tangent to an affine structure, $g$ is a flat Riemannian metric
(non--degenerate even symmetric quadratic tensor) such that
$\Cal{T}_M^f$ consists of $g$--flat tangent fields. Finally, $A$ is
an even symmetric tensor $A:\,S^3(\Cal{T}_M)\to \Cal{O}_M.$
All these data must satisfy the following conditions;

\smallskip

a) {\it Potentiality of $A$.} Everywhere locally there exists
a function $\Phi$ such that $A(X,Y,Z)=(XYZ)\Phi$ for
any flat vector fields $X,Y,Z.$

\smallskip

b) {\it Associativity.} $A$ and $g$ together define
a unique symmetric multiplication $\circ :\,\Cal{T}_M
\otimes \Cal{T}_M\to \Cal{T}_M$ such that 
$$
A(X,Y,Z)=g(X\circ Y,Z)=g(X,Y\circ Z).
$$
This multiplication must be associative.

\smallskip

In other words, in flat coordinates the tensor of
the third derivatives $\Phi_{ab}{}^c$ must constitute
the set of structure constants of an associative algebra.
 
\medskip

If one excludes the trivial case when $\Phi$ is a cubic form
with constant coefficients in flat coordinates,
the first large class of Frobenius manifolds
was discovered by Kyoji Saito even before
Dubrovin's axiomatization (see [S1], [S2] and [O]):

\smallskip
 
{\it (i) Moduli spaces of unfolding (germs of) isolated singularities
of hypersurfaces carry  natural structures of
Frobenius manifold.}

\smallskip

Each such structure is determined by a choice
of Saito's {\it good primitive form.}

\smallskip

In [D] a more global variation of this construction
is described (Hurwitz's spaces). 

\smallskip

Physicists call
the relevant TFT the topological sector of
the Landau--Ginzburg theory: cf. [C].

\medskip

The second large class of Frobenius manifolds
was discovered by physicists (Witten, Dijkgraaf, Vafa)
and is called Quantum Cohomology. For an axiomatic
treatment, see [KM1]. The correlators of this theory
are called Gromov--Witten (GW) invariants. 
Their actual construction in
the algebraic--geometric framework was carried out
in [B] following [BM] and [BF].

\smallskip

{\it (ii) The formal completion at zero of the cohomology
(super)space of any smooth projective or compact symplectic
manifold carries a natural structure of formal
Frobenius manifold.}

\medskip

The third large class of Frobenius manifolds was recently
constructed by S.~Barannikov and M.~Kontsevich ([BK]).

\smallskip

{\it (iii) The formal moduli spaces of solutions
to the Maurer--Cartan equations modulo gauge equivalence,
related to a class of the differential 
Gerstenhaber--Batalin--Vilkovisky (dGBV)
algebras, carry a natural structure of formal
Frobenius manifold.} 

\medskip

As their main application, Barannikov and Kontsevich construct
the dGBV--algebra starting with the Dolbeault complex
of an arbitrary Calabi--Yau manifold, and conjecture
that the resulting formal Frobenius manifold ($B$--model)
can be identified with the quantum cohomology
of the mirror dual Calabi--Yau manifold ($A$--model).
Another source of dGBV--algebras is provided by
the BRST cohomology of certain chiral algebras:
cf. [LZ].

\smallskip

The Mirror Conjecture is a part of the gradually
emerging considerably more general pattern.
Within the Calabi--Yau domain, it should be a
consequence of the Kontsevich's conjecture
about the equivalence of the Fukaya triangulated
category associated to one member
of the mirror pair and the derived category of sheaves
on the other member. Furthermore, one expects
the extension of the mirror picture
to other classes of varieties, non--necessarily smooth,
compact or having trivial canonical class.
For some exciting recent results on mirrors, cf. [Giv].

\smallskip

Isomorphisms of Frobenius manifolds of different classes
remain the most direct expression, although by no means the
final one, of various mirror phenomena.

\smallskip

From this vantage point, the three classes of 
examples considered above should be compared at least
in two ways.

\smallskip

First, one looks for isomorphisms between 
Frobenius manifolds (and their submanifolds)
constructed by different methods.

\smallskip

Second, one tries to generalize to other classes
of Frobenius manifolds additional structures
peculiar to each of the known classes. 

\medskip

Consider, for example, Quantum Cohomology.
Physically, quantum cohomology of a manifold $V$
is only the tree level small phase space part of the topological
sigma model with target space $V$. In particular, the
correlators of this theory, which are essentially
the coefficients of the formal Frobenius potential,
can be mathematically defined in terms of the
intersection indices on moduli spaces of stable maps
of curves of genus zero to $V$. This set of the correlators of Quantum
Cohomology of $V$ can be extended to a much vaster
structure involving, first, curves of arbitrary genus,
and second, the so called gravitational descendants,
mathematically expressible via Chern classes of certain
tautological bundles on the moduli spaces of stable maps
(cf. [KM2] for precise statements). 

\smallskip

This leads to two natural questions. 

\medskip

{\it What is the differential
geometric meaning of the arbitrary genus correlators with descendants
in the framework of general geometry of Frobenius
manifolds? 

\medskip

Which
Frobenius manifolds admit extensions of this type?}

\medskip

Of course, Frobenius manifolds  admitting such an extension
include those that are isomorphic to Quantum Cohomology.
Therefore a better understanding of this problem could shed some light
on the Mirror Conjecture as well. In [KM2] it is shown that
at least the correlators with gravitational descendants in any genus $g$
can be reconstructed
from the additional data consisting of two different parts.

\smallskip

One part is the genus $\le g$ Cohomological Field Theory
in the sense of [KM1] whose correlators take values
in the cohomology of the moduli spaces of stable curves
with marked points. At the moment it is unclear which abstract
Frobenius manifolds can be extended to such genus $\le g$
geometry and how it can be done. 

\smallskip

Another part of the data
concerns only genus zero correlators and therefore in principle
can be formulated in terms of arbitrary Frobenius manifold.
Its existence, however, poses  non--trivial restrictions 
on the manifold which are axiomatized below in the notion
of qc--type.

\medskip

This whole setup can be illuminated by comparison with the
motivic philosophy. In principle, any natural
structure on the cohomology of an algebraic manifold
can be considered as a realization of its motive,
and the question which abstract structures of a given type
arise from cohomology (``are motivated'') is a typical question
of the theory of motives.

\smallskip

Quantum Cohomology is a highly
nonlinear realization of the motives of smooth projective manifolds. 
It is functorial, at least in the
naive sense, only with respect to isomorphisms.
Hence it cannot be extended to the category of
the Grothendieck motives in an obvious way.
Nevertheless, the natural monoidal
structure of motives extends to Frobenius manifolds. 
Their tensor product in the formal context 
furnishing
the K\"unneth formula for Quantum Cohomology was
constructed in [KM] and [KMK]. 
R.~Kaufmann (cf. [K] and paper in preparation) has shown 
that the tensor product of convergent potentials converges,
and the resulting Frobenius manifold
in a sense does not depend on the choice of the base points.   
This adds some flexibility
to the motivic perspective. For example, Frobenius
manifolds provide
a  context in which one can meaningfully speak
about cohomology of fractional weight: cf. e.g. the treatment
of $A_n$--manifolds in 2.3.1 and 3.5 below.
It should be also
compared with  S.~Cecotti's suggestion that the
TFT's of Landau--Ginzburg type naturally give rise
to mixed Hodge structures (see [C]).

\medskip

{\bf 0.2. Plan of the paper.} This paper is dedicated to the study
of the Frobenius manifolds of the three classes (i), (ii) and (iii)
from the perspective described above.

\smallskip 

In \S 1 we start with reminding the
formalism of Gromov--Witten invariants and
the structure of the potential of the Quantum Cohomology.
In the treatment of [KM] its terms of degree $\le 2$
were not fixed because of absence of stable curves
of genus zero with $\le 3$ marked points.
However, the Divisor Axiom allows an unambiguous definition
of these terms (Proposition 1.3.1). This simple remark
is essential for the definition of $\le 2$ point
correlators for abstract formal Frobenius manifolds.
We introduce the notion of the manifold of qc--type,
which  embodies a version of the Divisor Axiom,
and show that it allows us to define for such manifolds
the correlators with gravitational descendants.
This is an elaboration of the picture sketched
in \S 3 of [KM2]. 

\medskip

In \S 2 we introduce the notion  of {\it the (strong) Saito framework.}
It axiomatizes those properties of the spaces
of miniversal deformations of isolated singularities
of functions which directly lead to the Frobenius
structure. The most difficult and deep
aspects of Saito's theory are thereby neatly
avoided and become ``existence theorems''.
(This illustrates the advantages of theft in comparison
with honest work, as was justly remarked 
about the axiomatic method in the beginning of this
century).

\smallskip

This part is taken from my notes to a lecture course
and is included here on suggestion of A.~Givental.

\smallskip

The main result of this section is summarized in the formulas
(2.12), (2.13) which in the context of Saito's
theory refer to the unfolding space of the direct sum
of singularities, and in our axiomatic treatment
are stated in terms of abstract direct sum diagrams.
Another proof of these formulas using
oscillating integrals was shown to me by A.~Givental.

\medskip

In \S 3 we show first of all, using (2.12) and (2.13), 
that the direct sum of Saito's frameworks
corresponds to the tensor
product of the associated Frobenius manifolds.

\smallskip

Looking then at the tensor products $M$ of the Frobenius
manifolds $A_n$ (deformation space of the singularity
$z^{n+1}$ at zero) and
more general manifolds with rational spectra we find out
that the integral part of their spectra define
Frobenius
submanifolds $HM$ which look like quantum cohomology
of a manifold with trivial canonical class,
at least on the level of discrete invariants (cf. below).
This argument exactly corresponds to the well known
idea of D.~Gepner ([Ge1], [Gr], [C]) of building Calabi--Yau
sigma models from the tensor products of minimal
models. The numerology is also similar.

\smallskip

Here is an example of our results.

\medskip

\proclaim{\quad 0.2.1. Claim} The manifold $A_n^{\otimes n+1}$
contains a (canonically defined) pointed Frobenius submanifold $HA_n^{\otimes n+1}$ whose spectrum looks formally like that of even--dimensional part of
quantum cohomology  of an $(n-1)$--dimensional
algebraic (or symplectic) manifold $V$.

\smallskip

More precisely, $V$ must have Betti numbers
$$
h^{2m}(V):=\roman{the\ number\ of\ } (i_1,\dots ,i_{n+1})\in\bold{Z}_{\ge 0}^{n+1}, 
\ \roman{satisfying}
$$
$$
\frac{1}{m}\,\sum_{k=1}^{n+1}i_k=m(n+1),\ 0\le i_k\le n-1
\eqno(0.1)
$$
and vanishing (modulo torsion) $c_1(V).$
\endproclaim

\smallskip

For example, even Betti numbers  must be (1,19,1) for $n=3$, and
(1,101,101,1) for $n=4$.

\smallskip

The Poincar\'e symmetry of them is generally established by the
involution $(i_k)\mapsto (n-1-i_k),\, m\mapsto n-1-m.$

\medskip

{\bf 0.2.2. Problem.} Is there actually a manifold $V_n$
whose Quantum Cohomology contains
$HA_n^{\otimes n+1}$?
Is it at least true that $HA_n^{\otimes n+1}$ is
Frobenius manifold of qc--type?
(Notice that $A_n$ itself is not of qc--type). 

\medskip

As was explained above, $A_n^{\otimes n+1}$ is the unfolding
space at zero of the singularity of $x_1^{n+1}+\ldots +x_{n+1}^{n+1}.$
 An argument which I learned from [C] 
and (in a different version) from A.~Givental then shows that
$HA_n^{\otimes n+1}$ carries the variation of  Hodge
structure corresponding to the middle cohomology
of the hypersurface $x_1^{n+1}+\ldots +x_{n+1}^{n+1}=0.$
More precisely, the volume form periods constitute the
horizontal sections of one of the structure connections
of the Frobenius manifold in question.

\medskip

Thus Problem 0.2.2 has the flavor of Mirror Conjecture, 
and of course Gepner's
idea was a precursor of the modern studies of the
mirrors. Hence at least the case $n=4$ of the Problem 0.2.2
must be reducible to the Givental's treatment
of the quintic fourfold [Giv].

\medskip

The Barannikov--Kontsevich construction [BK] conjecturally provides
another, and quite general, class of Frobenius manifolds
of Calabi--Yau $B$--type. To describe it succinctly,
notice that the space of vector fields on a Frobenius manifold
is simultaneously a Lie (super)algebra
and a (super)commutative algebra. The classical example of
such combination is the algebra of functions
on a Poisson (super)manifold. So it would appear that
Frobenius manifolds could come from Poisson
structures. However, such a relation cannot be straightforward,
because the Poisson identity between the bracket
and the multiplication does not hold on Frobenius 
manifolds. The ingenious twist in [BK]
consists in considering instead {\it odd
Poisson algebras  with two differentials}.
Under appropriate conditions, the Frobenius
manifold structure is then induced on the
homology space of such Poisson algebra.
 
\medskip
 
The sections 4--6 of this paper
constitute a completely self--contained account
of the theory of [BK] in the axiomatic context
of dGBV--algebras. Specifically, in \S 3 we supply direct
elementary proofs of all results related to the
formality and to the structure of Maurer--Cartan
moduli spaces at their non--obstructed points. 
In \S 4 we collect a list of basic  general properties of
dGBV--algebras. Finally, in \S 6 we define
the relevant formal Frobenius manifolds. We would like to stress
the similarity of the formulas defining
$\circ$--multiplication in the Saito's
and Barannikov--Kontsevich constructions:
compare (2.3) and (6.1). This supports the
expectation that both construction might be
special cases of a more general picture.

\medskip

As a comment to the title
of Cecotti's paper [C], it is instructive to compare the extensions of the Calabi--Yau
variations of Hodge structure (VHS) given by Saito's theory to the
construction of [BK] and the equivariant
theory of Givental ([Giv]).

\smallskip

Barannikov and Kontsevich embed any Calabi--Yau VHS into
a Frobenius supermanifold that has all discrete
invariants perfectly matching those of the mirror dual
quantum cohomology. Thus it has a good chance to be the correct
$B$--model in the classical Mirror Conjecture picture.
A drawback of this embedding is that
the relevant Frobenius manifolds are not semisimple,
and so the identification of them hardly can be
achieved by formal calculations.

\smallskip

To the contrary, whenever methods of Gepner, Saito and
Givental are applicable
(quasi--homogeneous
singularities, anticanonical hypersurfaces in toric
compactifications),  they embed (parts of) Calabi--Yau VHS
into Frobenius manifolds that are generically
semisimple or into families of manifolds
with generically semisimple general fiber. This makes more
accessible the direct check of mirror isomorphisms
(cf. \S 3 below).
But these techniques do not give the full dimension
spectrum (e.g. odd--dimensional cohomology is
skipped) and they are not directly applicable to
those Calabi--Yau manifolds which are not toric
anticanonical hypersurfaces (or complete intersections).

\medskip

{\it Acknowledgements.} The first draft of this paper
was written after stimulating discussions of
Gepner's work with V.~Schechtman. It was read and extensively 
commented by A.~Givental who corrected some mistakes and
supplied illuminating explanations of the Saito theory.
In particular, the whole \S 2 owes its existence to
Givental's suggestion. In writing the last three sections,
I have greatly benefited from M.~Kontsevich's
lectures at the MPI and his handwritten notes.  

\newpage

\centerline{\bf \S 1. Quantum cohomology and Frobenius manifolds.}

\bigskip 

{\bf 1.1. Gromov--Witten (GW) invariants.} We start with reminding
some basic notation and facts from [KM1], [BM], [B]. 

\smallskip

Let $V$ be a smooth projective algebraic manifold
over an algebraically closed field of characteristic zero,
$B=B(V)$ the semigroup of effective one--dimensional
algebraic cycles modulo numerical equivalence.
For any $\beta\in B(V),\,g,n\ge 0$ we can define
the Deligne--Mumford stack $\overline{M}_{g,n}(V,\beta )$
parametrizing stable maps of curves of genus $g$ with $n$ labelled
points, landing in class $\beta$. This stack comes equipped
with virtual fundamental class in the homological
Chow group with rational coefficients $J_{g,n}(V,\beta )\in 
A_s(\overline{M}_{0,n}(V,\beta ))$ where 
$$
s=(1-g)\,(\roman{dim}\,V-3)+(c_1(V),\beta )+n.
$$ 

\smallskip

Moreover, there are canonical morphisms $ev:\, \overline{M}_{g,n}(V,\beta )
\to V^n$ sending a stable map to the image of the family of
labelled points. In the stable range, that is
when $2g-2+n>0$, there is also
a map $st:\,\overline{M}_{g,n}(V,\beta )\to\overline{M}_{g,n}$
forgetting $V$ and stabilizing the curve.
They can be used in order to define {\it the Gromov--Witten
correspondences} in the Chow rings 
$$
I_{g,n}(V,\beta ):=(ev,st)_*(J_{g,n}(V,\beta ))\in
A_s(V^n\times\overline{M}_{g,n}).
$$
This family of Chow correspondences is the most manageable
embodiment of motivic quantum cohomology forgetting 
just the right amount of geometric information
encoded in the rather uncontrollable stack of stable maps.
For genus zero, the situation further simplifies. 
Since $V^n$ and $\overline{M}_{0,n}$
are smooth, we can identify $A_*$ with $A^*$. As S.~Keel
proved, $A^*(V^n\times\overline{M}_{0,n})=
A^*(V^n)\otimes A^*(\overline{M}_{0,n}),$
and $ A^*(\overline{M}_{0,n})$ is a finite--dimensional self--dual
linear space.
Hence one can identify $I_{0,n}(V,\beta )$
with the induced map $A^*(\overline{M}_{0,n})\to A^*(V^n).$
The space $A^*(\overline{M}_{0,n})$ is spanned by
the dual classes of the boundary strata $\overline{M}(\tau )$
indexed by $n$--trees. So calculating
$I_{0,n}(V,\beta )$ amounts to calculating a finite family
of elements in $A^*(V^n)$ indexed by these trees.

\medskip

{\bf 1.2. Frobenius manifolds.} All our examples of Frobenius manifolds 
(see 0.1) will come
equipped with two additional structures: a flat
vector field $e$ which is {\it identity} with respect
to $\circ$, and {\it an Euler vector field} $E$
expressing the scaling invariance of the Frobenius manifold $M$.
More precisely, we must have $\roman{Lie}_E(g )=Dg$
for some constant $D$, and $\roman{Lie}_E(\circ )=d_0\circ$
for another constant $d_0$ (in the context of Frobenius
manifolds $g$ means the metric, not the genus,
which is zero for the relevant GW--invariants). If the first condition is
satisfied, the second one is equivalent to
$$
E\Phi =(d_0+D)\Phi + \roman{a\ polynomial\ in\ flat\
coordinates\ of\ degree\ }\le 2.
$$

\smallskip

For any Euler field $E$ we have $[E,\Cal{T}_M^f]\subset
\Cal{T}_M^f.$ Assume for simplicity that the spectrum
$\{d_a\}$ of $-\roman{ad}\,E$ on flat vector fields
belongs to the base field. We understand $\{d_a\}$
as a family of constants with multiplicities.
The constant $d_0$ introduced earlier is 
in addition the eigenvalue corresponding to $e.$
The family $D,\{d_a\}$ is called the spectrum of $M.$
Since any multiple of $E$ is an Euler field together
with $E$, in the case $d_0\ne 0$ we can normalize $E$
by the condition $d_0=1$.

\smallskip

We also put $d=2-D$ and $q_a=1-d_a$ and call the
family $d,\{q_a\}$ {\it the $d$--spectrum} of $M$.
 
\smallskip

If $M$ is the formal spectrum of the ring of formal series
in flat coordinates, $\Phi$ is a formal solution of the Associativity
Equations, we call
$M$ {\it formal Frobenius manifold.} Formal Frobenius manifolds
can be tensor multiplied. The underlying metric space of flat fields
of the tensor product is the usual tensor product of the
respective spaces of factors. The potential of the tensor
product is defined in a much subtler way: see [KM1] and [M].
If the factors are additionally
endowed with flat identities  and normalized Euler fields,
they can be used to produce a canonical flat identity
and Euler field on the tensor product: see [K].

\medskip
 
{\bf 1.3. From genus zero GW invariants to Frobenius manifolds.}
In the situation of 1.1, we can construct a formal
Frobenius manifold $H^*_{\roman{quant}}(V)$ whose underlying linear supermanifold is
the completion of $H^*(V,\Lambda )$ at zero, with
obvious flat structure and Poincar\'e form as metric. Here
$\Lambda$ is a $\bold{Q}$--algebra endowed with the
universal character $B(V)\to\Lambda :\beta\mapsto q^\beta$,
with values in the Novikov ring $\Lambda$ which is the completed
semigroup ring of $B(V)$ eventually localized
with respect to the multiplicative system $q^{\beta}.$
It is topologically spanned by the  monomials
$q^{\beta}=q_1^{\beta_1}\dots\, q_m^{b_m}$ where
$\beta =(b_1,\dots ,b_m)$ in a basis of the numerical class
group of 1--cycles, and $(q_1,\dots ,q_m)$ are
independent formal variables. This is needed to make
convergent the formal sums
$$
I_{g,n}(V):=\sum_{\beta\in B}q^{\beta}I_{g,n}(V,\beta )\in A_{\Lambda}^*(V^n\times
\overline{M}_{g,n})
\eqno(1.1)
$$ 
For $H^*$ we can take any cohomology theory
functorial with respect to Chow correspondences.
The construction of the potential requires
only the top degree terms of the genus zero GW invariants. To be concrete,
choose a homogeneous basis $\{\Delta_a\,|\,
a=0,\dots r\}$ of $H^*(V,\bold{Q})$. Denote by
$\{x_{a}\}$ the dual coordinates
and by $\Gamma =\sum_{a}x_{a}\Delta_a$
the generic even element of the cohomology
superspace. This means that $\bold{Z}_2$--parity of
$x_a$ equals 0 (resp. 1), if $\Delta_a$ is even--dimensional
(resp. odd--dimensional). Put for $\gamma_i\in H^*(V)$
$$
\langle I_{g,n}(V)\rangle (\gamma_1\otimes\dots\otimes \gamma_n)=\langle \gamma_1\dots \gamma_n\rangle_{g,n} :=
\int_{I_{g,n}(V)} pr^*(\gamma_1\otimes\dots\otimes \gamma_n)
\eqno(1.2)
$$
where $pr:\, V^n\times\overline{M}_{g,n}\to V^n$ is the projection.
Then the quantum cohomology potential is
$$
\Phi (x)=\langle e^{\Gamma}\rangle_0:=
\sum_{n,(a_1,\dots ,a_n)}
\epsilon (a)\,\frac{x_{a_1}
\dots x_{a_n}}{n!}\,
\langle\Delta_{a_1}\dots\Delta_{a_n}
\rangle_{0,n}
\eqno(1.3)
$$
where $\epsilon (a)$ is the sign resulting from
rewriting $\prod x_i\Delta_{a_i}$ as $\epsilon (a)\prod x_i\prod \Delta_{a_i}.$
Assume that the dual fundamental class $\Delta_0$ of $V$
is the part of our basis. Then the flat identity
is $\partial_0=\partial /\partial x_0.$ Moreover, the
Euler field is
$$
E=\sum_a\left(1-\frac{|\Delta_a|}{2}\right)\,x^a\partial_a +
\sum_{b:\,|\Delta_b|=2}r^b\partial_b,
\eqno(1.4)
$$
where $\Delta_a\in H^{|\Delta_a|}(V)$, and $r^b$ are defined by
$$
c_1(\Cal{T}_V)=-K_V=\sum_{b:\,|\Delta_b|=2}r^b\Delta_b.
\eqno(1.5)
$$
Clearly, $-\roman{ad}\, E$ is semisimple on flat vector fields.
Let $H(d_a)$ be the eigenspace correspondng to $d_a$. We have $H(d_a)=H^{2q_a}(V).$
Hence the total spectrum is
$$
D=2-\roman{dim}\,V,\ d_0=1,\ d_a=1-\frac{|\Delta_a|}{2} 
\roman{\ of\ multiplicity\ dim\,}H^{|\Delta_a|}(V)
\eqno(1.6) 
$$
and the $d$--spectrum is
$$
d=\roman{dim}\,V,\ \{q\  \roman{\ of\ multiplicity\ } h^{2q}(V)\}.
\eqno(1.7)
$$
Quantum cohomology of the direct product 
of manifolds is the tensor product of the respective
Frobenius manifolds. So the $d$--spectrum behaves
as is expected.

\smallskip

$E$ itself is not flat, but it has the flat projection $E(0)$ to $H(D)=H^2$
which is just the anticanonical class, if $E$ is normalized as above
by $d_0=1$ and $H^*(V)$ is identified with the space of flat vector
fields. This is evident
from (1.3) and (1.4). 

\smallskip

The spaces $H(d_a),H(d_b)$ are orthogonal unless $d_a+d_b=D$, 
or equivalently  $q_a+q_b=d$, and dual
in the latter case. They also all have integral structure
compatible with metric. If we work with a coefficient ring $\Lambda$,
these subspaces are direct (free) sumbodules. Identity
belongs to $H(1)=H^0(V)$. Especially important are $H^2(V)=H(D)$,  $H_2(V)=H(0)$,
and the semigroup $B\subset H_{2,\bold{Z}}$ of effective algebraic
classes, in which every element
is finitely decomposable, and zero is
indecomposable. They are never trivial for projective smooth $V$
of positive dimension.

\smallskip

Returning to the potential $\Phi$, we see that since
$\overline{M}_{0,n}$ is empty for $n\le 2$, the definition (1.3)
specifies only its terms of degree $\ge 3$ in $x_a$. The
validity of the Associativity Equations is not
sensitive to this indeterminacy. However, the missing terms 
can be uniquely normalized either geometrically, by
integrating over $J_{0,n}(V,\beta )$ in $\overline{M}_{0,n}(V,\beta )$
for $n\le 2$, or formally, by using the Divisor Axiom
of [KM1]. Since this normalization is important
for the future use, we describe it explicitly. Denote by
$\delta :=\sum_{a:\,d_a=0}x_a\Delta_a$ the generic even element
of $H^2$.

\medskip

\proclaim{\quad 1.3.1. Proposition} There exists a unique formal
function $\Phi$ differing from (1.3) only by terms of
degree $\le 2$ which is representable as a formal Fourier
series in $q^{\beta}e^{(\beta ,\delta )},\,\beta \in B,$
with coefficients which are formal series of the
remaining coordinates, having the following properties.
Put $\Phi=\Psi +c$ where $c$ is the constant ($\beta =0$)
term of the Fourier series. Then, assuming $d_0=1$
and denoting by $E(0)$ the anticanonical class summand of $E$,
we have:

\smallskip

a) $E\Psi =(D+1)\Psi$.

\smallskip

b) $c$ is a cubic form with $(E-E(0))\,c=(D+1)c,$
the classical cubic self--intersection index divided by $6$.

\endproclaim

\medskip

In fact, if one puts formally
$q^{\beta}=0$ for $\beta \ne 0$ in the structure
constants of the quantum multiplication,
one gets the classical cup multiplication on $H^*(V)$ which together
with metric determines $c$ starting even with non necessarily
normalized potential.

\medskip

{\bf Proof.} We use the properties of $I_{0,n}(V,\beta )$ stated axiomatically
in [KM1] and proved in [B].

\smallskip

 The initial potential $\Phi$ is expressed
via Gromov--Witten invariants 
$$
\Phi(\Gamma ):=\sum_{n\ge 3}\sum_{\beta}
\,\frac{q^{\beta}}{n!}\langle I_{0,n}(V,\beta )\rangle
(\Gamma^{\otimes n})
$$
where $\langle I_{0,n}(V,\beta )\rangle$ is defined by the
same formula as (1.2) only with integral taken over
$I_{0,n}(V,\beta ).$ 
The part of $\Phi$
corresponding to $\beta =0$ is  exactly $c=\dfrac{1}{6}\,(\Gamma^3).$

\smallskip

The maps $\langle I_{0,n}(V,\beta )\rangle :\,H^*(V)^{\otimes n}\to
\Lambda$ are defined for $n\ge 3$ and satisfy
$$
\langle I_{0,n}(V,\beta )\rangle (\alpha\otimes\delta )=
(\beta ,\delta )\,\langle I_{0,n-1}(V,\beta )\rangle (\alpha )
$$
for $\delta \in H^2(V)$ (we write $(\beta ,\delta )$
for $g(\beta ,\delta )).$ This follows from the Divisor Axiom.
It is easy to check that there exists
a unique polylinear extension of $\langle I_{0,n}(V,\beta )\rangle$ to
all $n\ge 0$ satisfying this identity. In fact, it suffices to
put
$$
\langle I_{0,n}(V,\beta )\rangle (\alpha )=
(\beta ,\delta )^{-m} \langle I_{0,n+m}(V,\beta )\rangle (\alpha
\otimes \delta^{\otimes m})
$$
for any $(m,\delta )$ with $m+n\ge 3$ and invertible
$(\beta ,\delta )$.

\smallskip

Now put
$$
\Psi(\Gamma ):=\sum_{n\ge 0}\sum_{\beta \ne 0}
\,\frac{q^{\beta}}{n!}\langle I_{0,n}(V,\beta )\rangle
(\Gamma^{\otimes n})
$$
Clearly, $\Psi +c$ differs from the initial $\Phi$ by terms of degree
$\le 2.$ Moreover, $c$ is a linear combination of
 $x^ax^bx^c$
with $|\Delta_a|+|\Delta_b|+|\Delta_c|=2\,\roman{dim}\,V$ so that $(E-E(0))\,c=(D+1)c.$
As for $E\Psi$, we have for $\Gamma =\gamma_0 +\delta$
$$
\Psi (\gamma_0+\delta )=
\sum_{i,k\ge 0}\sum_{\beta \ne 0}
\,\frac{q^{\beta}}{i!\,k!}\langle I_{0,n}(V,\beta )\rangle
(\gamma_0^{\otimes i}\otimes\delta^{\otimes k})=
$$
$$
\sum_{i\ge 0}\sum_{\beta \ne 0}
\,\frac{q^{\beta}e^{(\beta ,\delta )}}{i!}\langle I_{0,n}(V,\beta )\rangle
(\gamma_0^{\otimes i}).
\eqno(1.8)
$$
Let us apply now $E$ to any summand in (1.8). The $E(0)$ part
acts only upon $e^{(\beta ,\delta )}$ and multiplies it
by $(c_1(V),\beta ).$ The $E-E(0)$ part multiplies any
monomial $x_{a_1}\dots x_{a_n}$ in non--divisorial
coordinates by $\sum_i(1-|\Delta_{a_i}|/2).$ From (1.2)
and (1.3) we see that $\beta$ can furnish a non--zero contribution
to such term only if
$$
\roman{dim}\,I_{g,n}(V,\beta )=\roman{dim}\,V-3+(c_1(V),\beta )+n=
\sum_{i=1}^n \frac{|\Delta_{a_i}|}{2}.
$$
Hence every non--vanishing term of (1.8) is an eigenvector
of $E$ with eigenvalue $D+1=3-\roman{dim}\,V.$

\smallskip

This proves the Proposition.

\smallskip

Notice in conclusion that $q^{\beta}e^{(\beta ,\delta )}$
is the universal character of $B$ together with $q^{\beta}.$
We have introduced $q^{\beta}$ only to achieve the formal convergence.
If it holds without $q^{\beta}$, we can forget about it.
Moreover, if the formal Fourier series actually converges
for $\delta$ lying somewhere in the complexified
ample cone, $\Psi (x)$ has a free abelian symmetry
group: translations by an appropriate discrete subgroup
in the space $H^2(V,i\bold{R}).$ Conversely, in the analytic
category this condition is necessary for
the existence of the appropriate Fourier series. 

\medskip

{\bf 1.4. Potentials of qc--type.} Based upon the analysis above,
we will introduce the following definition.  Its first goal is
to axiomatize a part of the structures of sec. 1.3
which suffices for the construction of 
the coupling of a formal Frobenius manifold with
gravity in the sense of [KoM2].
As we will recall below, this construction
is divided into two steps: the construction of the
modified gravitational descendants which can be done
for any formal Frobenius manifold, and the construction
of a linear operator $T$ on the big phase space 
which requires additional structures.

\smallskip

The second goal is to provide an intermediate step
in the problem of checking whether a given
formal Frobenius manifold is quantum cohomology.
We must be able at least to detect the following
structures. 

\medskip

\proclaim{\quad 1.4.1. Definition} Let $(M=\roman{Spf}\,k[[H^t]],g ,\Phi_0 )$
be a formal Frobenius manifold over a $\bold{Q}$--algebra $k$
with flat identity, Euler field $E$, 
and spectrum $D,\{d_a\}$ in $k$ as above. 
Here $H$ is a free $\bold{Z}_2$--graded $k$--module
of flat vector fields, and $H^t$ is the dual
module of flat coordinates vanishing at the origin.
Put $H^2=H(0),H_2=H(D)$.

\smallskip

Assume that there exists a semigroup $B\subset H_2$
with finite decomposition and indecomposable zero, and the cubic form
$c$ on $H$, such that by eventually changing terms of degree $\le 2$
in $\Phi_0$ we can obtain the potential of the form 
$$
\Phi=\Psi +c,\,E\Psi=(D+1)\Psi,\,(E-E(0))c=(D+1)c,
$$
$$
\Psi (\gamma_0+\delta )= \sum_{i\ge 0}\sum_{\beta\in B\setminus\{0\}}
\frac{e^{(\beta ,\delta )}}{i!}\,I_{\beta}(\gamma_0^{\otimes i})
\eqno(1.9)
$$
such that all summands in the last sum are eigenvectors of
$E$ with eigenvalue $D+1$. Here $\gamma$ is a
generic even element of $H$, $\delta$ its ``divisorial''
$H^2$--part, $\gamma_0=\gamma -\delta .$ The coefficient
$I_{\beta}(\gamma_0^{\otimes i})$ is a form in non--divisorial
coordinates.

\smallskip

A formal Frobenius manifold satisfying these conditions will
be called  of qc--type.

\endproclaim

\smallskip

A flat identity $e$ in this language is an element $e\in H$
which considered as a derivation satisfies
$$
e\,\Psi=0,\,e\,c=g
$$
(see [M], p. 29 for the same expressions in coordinates).

\medskip

{\bf 1.4.2. Correlators of qc--manifolds.} Let $M$ be a formal Frobenius
manifold of qc--type. Recall that $\Phi_{ab}{}^c$ are the structure
constants of the quantum multiplicarion.
On qc--manifolds there are two useful specializations
of this structure.

\smallskip

a) The ``small quantum multiplication''  obtained by 
restricting $\Phi_{ab}{}^c$ to $\gamma_0=0$.
 We will denote this multiplication by dot.

\smallskip

b) The cup multiplication $\cup$ obtained by putting formally
$e^{(\beta, \gamma )} =0$ for all $\beta\ne 0$ (``large volume limit'').
In other words, this is the multiplication, for which $c$
can be written as
$$
c(\gamma )=\frac{1}{6}\, g(\gamma , \gamma\cup\gamma ).
$$
We now define correlators $\langle\dots\rangle :\,H^{\otimes n}\to k$ as $S_n$--invariant
polylinear
functions whose values  are derivatives of $\Phi$
at zero. In other words, for a basis $\{\Delta_a\}$ of $H$
and dual coordinates $\{x_a\}$ as above, we have
$$
\Phi (x)=\sum_{n,a_1,\dots ,a_n}
\epsilon (a)\,\frac{x_{a_1}
\dots x_{a_n}}{n!}\,
\langle \Delta_{a_1}\dots \Delta_{a_n}\rangle  
$$
In the qc--case we can write
$$
\langle \Delta_{a_1}\dots \Delta_{a_n}\rangle =
\sum_{\beta\in B\setminus\{0\}}\langle \Delta_{a_1}\dots \Delta_{a_n}\rangle_{\beta} + \langle \Delta_{a_1}\dots \Delta_{a_n}\rangle_0
$$
where the first sum comes from $\Psi$ and the second, nonvanishing
only for triple arguments, from $c$.

\smallskip

Looking at (1.9) one sees, that small quantum multiplication
depends only on the triple correlators of non--divisorial
elements of the basis.

\medskip

\proclaim{\quad 1.4.3. Claim} The correlators of the Frobenius manifolds
of qc--type satisfy the following Divisor Identity: if $\delta\in H^2$,
$\beta\ne 0,$
$$
\langle\delta \gamma_1\dots\gamma_n\rangle_{\beta}=
(\delta ,\beta )\langle\gamma_1\dots\gamma_n\rangle_{\beta}.
$$
\endproclaim
\smallskip

Reading backwards the proof of (1.8), one sees that this property
follows from (1.9).

\smallskip

This formula allows us to extend the definition of the correlators
to $n\le 2$ arguments.

\medskip

{\bf 1.5. Gravitational descendants for the Frobenius manifolds of
qc--type.} Let now $M$ be a formal Frobenius manifold as in 1.4.1
whose space we identify with the linear superspace $H$,
At first we do not assume that it is of qc--type.
Following [KM2], we can define its  modified
correlators with gravitational descendants. They are
polylinear functions on {\it the big phase space}
$\oplus_{d\ge 0}H[d]$ where $H[d]$ are copies of the
space $H$ identified with the help of the
shift operator $\tau :\,H[d]\to H[d+1].$ To define
them explicitly, we recall that any formal
Frobenius manifold gives rise to the genus zero
Cohomological Field Theory.
Namely, there exists a unique
sequence of linear maps
$I^M_n:\,H^{\otimes n}\to H^*(\overline{M}_{0,n},k),\,n\ge 3,$
satisfying the folowing properties.

\smallskip

a) $I_n^M$ are $S_n$--invariant and compatible with restriction to the
boundary divisors (cf. [KM1] or [M], p. 101).

\smallskip

b) The top degree term of $I^M_n$ capped with the fundamental 
class of $\overline{M}_{0,n}$ is the correlator of $M$ with $n$ arguments.

\smallskip

Moreover, in the quantum cohomology case
$$
I^M_n=\sum_{\beta}q^{\beta}I^V_{0,n,\beta}
$$
where $I^V_{0,n,\beta}$ are the genus zero Gromov--Witten
invariants.

\smallskip
 
Now let $C\to \overline{M}_{0,n}$ be the universal curve,
$s_i:\,\overline{M}_{0,n}\to C,\,i=1,\dots ,n$ its
structure sections, $\omega_C$ the relative dualizing sheaf,
$\psi_i:=c_1(s_i^*(\omega_c))\in A^1(\overline{M}_{0,n}.$
The {\it the modified correlators with gravitational
descendants} for $M$ are defined by the formula
$$
\langle \tau^{d_1}\Delta_{a_1}\dots \tau^{d_n}\Delta_{a_n}\rangle :=
\int_{\overline{M}_{0,n}}I_n^M(\Delta_{a_1}\otimes\dots 
\otimes\Delta_{a_n})\psi_1^{d_1}\dots\psi_n^{d_n}
\eqno(1.10)
$$
and the generator function for them, {\it the modified
potential}, by the formula
$$
G^M(x)
= \sum_{n\ge 3,(a_i,d_i)}
\epsilon (a)\,\frac{x_{d_1,a_1}
\dots x_{d_n,a_n}}{n!}\,
\langle \tau^{d_1}\Delta_{a_1}\dots\tau^{d_n}\Delta_{a_n}
\rangle 
\eqno(1.11)
$$
where $(x_{d,a})$ are coordinates dual to $\tau^d\Delta_a.$

\smallskip

If we assume in addition that $M$ is of qc--type,
we can define also the (unmodified) two argument correlators
$\langle\tau_d \gamma_1\,\gamma_2\rangle$ by the inductive formula (25)
from [KM2]:
$$
\langle\tau_d\gamma_1\,\,\gamma_2\rangle =
\sum_{j=1}^{d}(-1)^{j+1}
\partial_{\delta}^{-j}
\langle\tau_{d-j}\gamma_1\,\,\delta\cdot
(\delta^{j-1}\,\cup\,\gamma_2)\rangle +
$$
$$
(-1)^d
\partial_{\gamma_0}^{-(d+1)}[
\langle\delta\,\,\gamma_1\,\,
\delta^d\,\cup\,\gamma_2\rangle
-\langle\delta\,\,\gamma_1\,\,
\delta^{d}\cup\gamma_2\rangle_{0}].
\eqno(1.12)
$$
Here $\delta\in H^2$ is an arbitrary (say, generic) element
such that $(\delta ,\beta )\ne 0$ for all $\beta\in B\setminus\{0\}$
and the operator $\partial_{\delta}^{-1}$ divides
$\langle\dots\rangle_{\beta}$ by $(\delta ,\beta ).$

\smallskip

Furthermore, put
$$
y_{c,b}=x_{c,b}+\sum_{(a,d),g\ge c+1}x_{d,a}\langle
\tau_{d-c-1}\Delta_a\,\Delta^b\rangle .
$$
Then the big phase space potential of $M$ is,
by definition, $F^M(x):=G^M(y)$, and the unmodified
correlators with gravitational descendants
of $M$ are defined as coefficients of $F$:
$$
F^M(x)
= \sum_{n\ge 3,(a_i,d_i)}
\epsilon (a)\,\frac{x_{d_1,a_1}
\dots x_{d_n,a_n}}{n!}\,
\langle \tau_{d_1}\Delta_{a_1}\dots\tau_{d_n}\Delta_{a_n}
\rangle 
\eqno(1.13)
$$
The main result of [KM2] is that if $M$ is  the quantum
cohomology of $V$, this prescription provides
the correlators with descendants of the topological sigma
model with target space $V.$ The latter are defined
by the formula similar to (1.10) but with
$\overline{M}_{0,n}$ replaced by $\overline{M}_{0,n}(V),$
$I$ replaced by $J,$ and the respective change in the meaning
of $\psi_i.$

\smallskip

In conclusion  notice that the sigma model
correlators satisfy, partly demonstrably,
partly conjecturally, some additional identities,
of which the most interesting are probably
the Virasoro constraints. I do not know
which of these identities might be valid
for the more general qc--type manifolds.
In any case, it would be interesting to
determine differential equations at least
for the modified potential with descendants
for general formal Frobenius manifolds or
particular examples like the manifolds  $A_n$
(see 2.3 below).

\newpage

\centerline{\bf \S 2. K.~Saito's frameworks}

\bigskip

{\bf 2.1. Setup.} Let $p :\,N\to M$ be a submersion of
complex analytic or algebraic manifolds, 
generally non--compact, $F$ a holomorphic function
on $N.$ We consider $F$ as a family of functions on the fibers
of $p$ parametrized by points of $M.$ In local
coordinates $z=(z_a), t=(t_b)$ where $t_b$ are constant
along the fibers of $p$ we write $F=F(z,t).$
 
\smallskip

Let $d_p:\,\Cal{O}_N\to\Omega^1_{N/M}$ be the relative
differential. Denote by $C$ the closed analytic subspace
(or subscheme) of the critical points of the restrictions of $F$
to the fibers given by the equation $d_pF=0.$ Its
ideal $J_F$ is locally generated by the partial derivatives
$XF$ where $X$ are vertical vector fields on $N.$
Derivatives $\partial F/\partial z_a$ of course suffice.
Let $i_C:\,C\to N$ be the natural embedding,
$p_C$ the restriction of $p$
to $C.$ 

\smallskip

Denote by $\Omega^{max}_{N/M}$  the invertible sheaf of holomorphic
vertical volume forms on $N,$ $L:=i_C^*(\Omega^{max}_{N/M}).$
The Hessian $\roman{Hess}\,(F)\in \Gamma (C,L^2)$ is a well defined section of
$L^2$ which in local coordinates as above can be written as
$$
\roman{Hess}\,(F)=i_C^*\,\left[ \roman{det}\,\left(\frac{\partial^2F}{\partial z_a\partial z_b}\right)\,(dz_1\wedge \dots\wedge dz_n)^2\right].
$$
We denote by $G_C\subset C$ the subspace $\roman{Hess}\,(F)=0.$
Let $\Cal{T}_M$ be the tangent sheaf of $M.$
Finally, let $\omega$ be a nowhere vanishing global section of
$\Omega^{max}_{N/M}.$

\medskip

\proclaim{\quad 2.1.1. Definition} The family
$(p :\,N\to M, F, \omega )$ is called  Saito's
framework, if the following conditions are satisfied:

\smallskip

a) Let the map $s:\, \Cal{T}_M\to p_{C*}(\Cal{O}_C)$ be
defined by ${X}\mapsto \overline{X}F\, \roman{mod}\,J_F,$
where $\overline{X}$ is any local (in $N$) lift of $X.$
Then $s$ is an isomorphism of $\Cal{O}_M$--modules.
In particular, $C$ is finite and flat over $M.$
Assume moreover that $G_C$ is a divisor, and $p_C:\,C\to M$ is \'etale
outside the complement to the divisor $G=i_{C*}(G_C)\subset M.$

\smallskip

b) Define the following 1--form $\epsilon$ on $M\setminus G.$
Its value on the vector field $X=s^{-1}(f)$ corresponding to the
local section  $f$ of $\pi_*(\Cal{O}_C)$ equals
$$
i_X(\epsilon ): = \roman{Tr}_{C/M}\,\left(\frac{f\,i_C^*(\omega^2)}{ \roman{Hess}\,(F)}\right)=
\sum_{i=1}^{\mu}\frac{f(\rho_i)}
{\roman{det}\,((\partial^2F/\partial z_a \partial z_b)(\rho_i))}
\eqno(2.1)
$$
where $\rho_i$ are the local branches
of the critical locus $C$ over $M,$ $(z_a)$
is any vertical local coordinate system unimodular with respect
to $\omega .$

\smallskip

Then the scalar product $g:\,S^2(\Cal{T}_{M\setminus G})\to \Cal{O}_{M\setminus G}$ defined by
$$
g(X,Y):=
i_{X\circ Y}(\epsilon )
\eqno(2.2)
$$
is a flat metric. Both $\epsilon$ and
$g$  (as flat metric) extend regularly to $M.$

\endproclaim

\medskip

{\bf 2.2. The (pre--)Frobenius structure associated to
the Saito framework.} Let $(p :\,N\to M, F, \omega )$
be a Saito framework.

\smallskip

Denote by $\circ$ the multiplication in $\Cal{T}_M$
induced by the one in $p_*(\Cal{O}_C)$:
$$
X\circ Y:=s^{-1}(\overline{X}F\cdot \overline{Y}F)\,\roman{mod}\,J_F
$$
or equivalently
$$
\overline{X\circ Y} F \equiv
\overline{X}F\cdot \overline{Y}F\,\roman{mod}\,J_F .
\eqno(2.3)
$$
Clearly, the vector field $e:= s^{-1}(1\,\roman{mod}\,J_F)$
is the identity for $\circ$. 
Let $\Cal{T}_M^f$ be the sheaf of vector fields
flat with respect to $g.$ Finally, put
$$
A(X,Y,Z)=g(X\circ Y,Z)=g(X,Y\circ Z).
$$
The last equality follows from (2.2) and the associativity
of $\circ .$ The tensor $A$ is symmetric because
$\circ$ is commutative. Therefore we have:

\medskip

\proclaim{\quad 2.2.1. Claim} The data $(\Cal{T}^f_M,g,A)$
define on $M$ the structure of pre--Frobenius manifold
in the sense of [M], Ch. I, Def. 1.1.1.
\endproclaim

\medskip

One can say more about this structure restricted to
$M\setminus G.$ Call a connected open subset $U$ in $M\setminus G$
{\it small} if $p_C^{-1}(U)$ is the disjoint union of
$\mu =\roman{dim}\,M$ connected components $U_i$ canonically isomorphic to $U.$
For concreteness, we will arbitrarily number
them by $\{1,\dots ,\mu\}$ as in (2.1). Then we have natural
ring isomorphisms
$$
\Gamma (U,p_{C*}(\Cal{O}_C))=\oplus_{i=1}^{\mu}
\Gamma (U_i,\Cal{O}_C)=\Gamma (U,\Cal{O}_M)^{\mu}.
$$
This $\Gamma (U,\Cal{O}_M)$--algebra has a basis of idempotents $f_i:=\delta_{ij}$ on $U_j.$
Defining $e_i\in \Gamma (U,\Cal{T}_M)$ by $s(e_i)=f_i,$
that is, $\overline{e}_iF\,\roman{mod}\,J_F=f_i,$
we get a local $\Cal{O}_M$--basis of $\Cal{T}_M$
satisfying $e_i\circ e_j=\delta_{ij}e_j$ and $e=\sum_i e_i.$
Denote by $u^i\in \Gamma (U,\Cal{O}_M)$ the restriction
of $F$ to $U_i$ pushed down to $U$ that is, put 
$$
u^i=F(\rho_i)
$$
in the notation of (2.1).
Small subsets cover $M\setminus G$ so that the structure group
of $\Cal{T}_M$ is reduced to $S_{\mu}.$ Summarizing, we have:

\medskip    
 
\proclaim{\quad 2.2.2. Proposition} The data $(\Cal{T}^f_M,g,A)$
define on $M\setminus G$ the structure of semisimple pre--Frobenius manifold
in the sense of [M], Ch. I, Def. 3.1 and 3.2.
Moreover, we have $e_iu^j=\delta_{ij}$ so that
$(u^j)$ form a local coordinate system (Dubrovin's canonical
coordinates) and $[e_i,e_j]=0$ because $e_i=\partial /\partial u^i.$
\endproclaim

\medskip

{\bf Proof.} Only the last statement might need
some argumentation. We have $p^*(e_iu^j)=\overline{e}_ip^*(u^j)$
for any lift $\overline{e}_i$ of $e_i.$ To calculate the
right hand side we can restrict it to any local section of
$p$ since it is constant along the fibers. We choose
$\overline{e}_i$ tangent to $U_j$ and restrict the right
hand side to $U_j$ where $p^*(u^j)$ coincides with $F.$
The result is $\delta_{ij}$ by the definition of $e_i.$

\medskip

For the future use, we can reformulate this as follows.
Dualizing $s$ we get the isomorphism 
$s^t:\,\Omega^1_M\to\Cal{H}om_{\Cal{O}_M}\,
(p_{C*}(\Cal{O}_C),\Cal{O}_M).$ Then $s^t(du^i):\,
p_{C*}(\Cal{O}_C)\to\Cal{O}_M$ is
the map which annihilates $j$--components for $j\ne i$
and coincides with the pushforward on the $i$--th component. 

\medskip

\proclaim{\quad 2.2.3. Theorem} The structure $(M,\Cal{T}^f_M,g,A)$
associated to the Saito framework is Frobenius iff $d\epsilon = 0.$
\endproclaim

\smallskip

{\bf Proof.} To check the Frobenius property on $M\setminus G$
we appeal to [M], Ch. I, Th. 3.3 (Dubrovin's criterium),
both conditions of which,
$[e_i,e_j]=0$ and $d\epsilon=0$ are satisfied. To pass from $M\setminus G$
to $M$ one can use a continuity argument, e. g. in the
following form, again due to Dubrovin.  Let $\nabla_0$
be the Levi--Civita connection of $g$, and
$\nabla_{\lambda}$ the pencil of connections on $\Cal{T}_M$
determined by its covariant derivatives 
$\nabla_{\lambda ,X}(Y):=\nabla_{0,X}(Y)+X\circ Y.$
Then $M$ if Frobenius iff $\nabla_{\lambda}$ is flat for some $\lambda \ne 0,$
and so automatically for all $\lambda .$
Clearly, this is the closed property.

\medskip

We will now discuss when $e$ is flat.

\smallskip

On a small $U$, we can define functions $\eta_j$ by 
$\eta_j =i_{e_j}(\epsilon )=g(e_j,e_j).$ When $e_j=\partial /\partial u^j$,
the closedness of $\epsilon =\sum_i\eta_idu^i$ means that $\eta_j=e_j\eta$
for a local function $\eta$ well defined up to addition
of a constant, or else $\epsilon =d\eta .$ In the notations (2.1)
$$
\eta_i =\frac{1}{\roman{det}\,((\partial^2F/\partial z_a \partial z_b)(\rho_i))}.
\eqno(2.4)
$$

\medskip

\proclaim{\quad 2.2.4. Theorem} Assume that the conditions of the 
Theorem 2.2.3 are satisfied.

\smallskip

The identity $e$ is flat, iff for all $i$, $e\eta_i =0,$
or equivalently, $e\eta = g(e,e) = \roman{const}.$
This holds automatically in the presence
of an Euler field $E$ with $D\ne 2d_0$ (see 1.2 above and 2.2.5 below).
\endproclaim

This is Prop. 3.5 from [M], Ch. I.

\medskip

One  important remark about the identity is in order.
Namely, in all examples I know of there exists a lift
$\overline{e}$ of $e$ to $N$ such that $\overline{e}F=1$
identically, so that in the appropriate coordinate system we have
$F=F_0+t_{0}$ where $F_0$ does not depend on $t_0$
which is lifted from $M$, and $\overline{e}=\partial /\partial t_0.$

\medskip

It remains to clarify what Euler fields  this structure can have.

\medskip

\proclaim{\quad 2.2.5. Theorem} Assume that the conditions of the 
Theorem 2.2.3 hold. Let $E$ be a vector field
on a small subset $U$ in $M.$

\smallskip

a) We have $\roman{Lie}_E(\circ )=d_0\circ$ iff 
$$
E=d_0\sum_i(u^i+c^i)e_i
\eqno(2.5)
$$
for some constants $c^i,$ where $(du^i)$ are 1--forms dual to $(e_j).$

\smallskip

In particular, for non--zero $E$
we have $d_0\ne 0$ so that we may normalize $E$ by $d_0=1$.
Furthermore, if the monodromy representation of the fundamental group
of $M\setminus G$ on $H_0$ of the fibers of
$C\to M$ has only one-dimensional trivial subrepresentation, the global vector
field $E$ of this form with fixed $d_0$ is defined uniquely
up to addition of a multiple of $e$.

\smallskip

b) For a field $E$ of the form (2.5) and a constant $D$,
we have $\roman{Lie}_E(g)=Dg$ iff
$$
E\eta =(D-d_0)\eta + \roman{const}.
\eqno(2.6)
$$
In particular, if $e$ is flat, adding a multiple of $e$ does not change
the validity of this property.
\endproclaim

\medskip

This follows from [M], Ch. I, Th. 3.6. 

\medskip

When $M$ comes from the Saito framework, we have a natural
candidate for the global Euler field with $d_0=1$ suggested by our
identification of local coordinates $u^i$.
Namely, put on any small $U$
$$
E_F:= \sum_{i=1}^{\mu}F(\rho_i)e_i = \sum_{i=1}^{\mu}u^ie_i.
\eqno(2.7)
$$
Assume that it is in fact an Euler field and that we are in the conditions
when it is defined uniquely up to a shift by
a multiple of $e$. Assume furthermore that there exists
a point $0$ in $M$ to which $E_F$ extends and at which
it vanishes ($0$ may lie in
$G,$ and in the theory of singularities it does so).
Since $e$ cannot vanish, the choice of such $0$
fixes $E_F$ completely.

\medskip

\proclaim{\quad 2.2.6. Definition} Saito's framework
$(p :\,N\to M, F, \omega )$ is called the strong Saito
framework, if the structure $(M, \Cal{T}^f_M, g, A)$
described above
is Frobenius, with flat identity $e$ and Euler field $E_F.$
\endproclaim

\medskip

{\bf 2.2.7. Remark.} Since the definitions of the pre--Frobenius
and Frobenius structures, and also of the identity and Euler fields,
are local, we can lift all these structures from
$M\setminus G$ to $C\setminus G_C.$

\medskip

{\bf 2.3. Unfolding singularities.} K. Saito's theory
(cf. [S1], [S2], [O] and the references therein)
produces (a germ of) a strong Saito's framework starting with
a germ of holomorphic function $f(z_1,\dots ,z_n)$ with isolated
singularity at zero.

\smallskip

Namely, one can choose holomorphic germs $\phi_0=1,\phi_1,
\dots ,\phi_{\mu -1}$ whose classes constitute a basis of
the Milnor ring $\bold{C}\{\{z\}\}/(\partial f/\partial z_a)$ 
in such a way that $F:=f+\sum t_i\phi_i$ is the miniversal unfolding of 
$f$. Then $N=N_f,$ resp. $M=M_f$ is a neighborhood of zero 
in the $(z,t)$--, resp. $(t)$--space, and $F$ is defined above.

\smallskip

The crucial piece of the structure is the choice
of $\omega$ encoded in the Saito notion of
a good primitive form. Generally its existence
is established in an indirect way. However, if
$f$ is a quasi--homogeneous polynomial, one can take
$\omega =dz_1\wedge \dots\wedge dz_n.$
In this case most of the data constituting the
Saito framework are algebraic
varieties, rational maps and rational
differential forms so that the whole setup
has considerably more global character.

\smallskip

In order to help
the reader to compare notation, we notice that
Saito's $S$ is our $M$, and 
our function $\eta$ is denoted
$\tau$ on p. 630 of [O]. Starting with the
germ of $z^{n+1}$ at zero, one obtains in this way
the following Frobenius manifold.

\medskip

{\bf 2.3.1. Example: manifolds $A_n$.} Denote by
$N$ , resp. $M$, the affine space with coordinates
$(z;a_1,\dots ,a_n)$, resp. $a_1,\dots ,a_n)$, and
by $p$ the obvious projection. Put
$$
F=F(z;a_1,\dots ,a_n)=z^{n+1}+a_1z^{n-1}+\dots +a_n.
$$
Then $C$ iz given by the equation $F^{\prime}(z)=0.$
We choose $\omega =d_pz.$

\smallskip

Making explicit the basic structures described above
we get the following description of $A_n.$
Consider the global covering of $M$ whose points
consist of total orderings of the roots $\rho_1,\dots ,
\rho_n$ of $F^{\prime}(z).$ On
the semisimple
part of it where $F^{\prime}(z)$ has no multiple
roots and $u^i:=F(\rho_i)$ are local coordinates
we have the flat metric
$$
g:=\sum_{i=1}^n\frac{(du^i)^2}{F^{\prime\prime}(\rho_i)}
$$
with metric potential
$$
\eta =\frac{a_1}{n+1}=\frac{1}{n-1}\sum_{i<j}\rho_i\rho_j=
-\frac{1}{2(n-1)}\sum\rho_i^2.
$$
\smallskip

Furthermore, $e,\,E$ and flat coordinates $x_1,\dots ,x_n$
can be calculated through $(a_1,\dots ,a_n)$ (which are
generically local coordinates as well):
$$
e=\partial /\partial a_n,\ \roman{i.~e.,\ } ea_n=1, ea_i=0 
\roman{\ for\ } i<n.
$$
$$
E=\frac{1}{n+1}\sum_{i=1}^{n}(i+1)a_i\frac{\partial}{\partial a_i},
$$
$x_i$ are the first Laurent coefficients of the inversion of
$w={\root n+1 \of{p(z)}} = z+O(1/z)$ near $z=\infty :$
$$
z=w +\frac{x_1}{w} + \frac{x_2}{w^2} +\dots +\frac{x_n}{w^n}
+ O(w^{-n-1}).
$$

\smallskip

For the direct proof of these statements, see e.g. [M], Ch. I, 4.5.

\smallskip

The spectrum of $A_n$ is $D=\dfrac{n+3}{n+1},\,
d^{(i)}=\dfrac{i+1}{n+1},\, 1\le i\le n,$ more precisely,
$Ex_i= \dfrac{i+1}{n+1}x_i.$ Hence the $d$--spectrum is
$$
d^{(n)}=\frac{n-1}{n+1},\, q_i^{(n)}=\frac{i}{n+1},\,i=0,\dots ,n-1.
$$
Now, $\Phi$ is analytic in $x_a$ and the spectrum of $-\roman{ad}\,E$
is strictly positive. Therefore $\Phi$ must be a polynomial
in flat coordinates. One can check that its degree
is precisely $n+2$. Hence for $n\ge 2$ it cannot be
of qc--type, and by the  method described in \S 1 we can define for $A_n$
only modified correlators with gravitational descendants.

\smallskip

Comparing the spectrum of $A_n$ with that of the quantum cohomology
of projective spaces, one can somewhat imaginatively
say that $A_n$ represents ``projective space of dimension
$\dfrac{n-1}{n+1}$, with rank one cohomology
in each dimension $\dfrac{i}{n+1}, \,0\le i\le n-1$
and with vanishing canonical class''.

\medskip

{\bf 2.3.2. Example: Gepner's manifolds $V_{n,k}$.} Let $n\ge 2,k\ge 1,
h=n+k.$ We will call {\it Gepner's Frobenius manifold}
$V_{n,k}$ the manifold which is produced from the Saito's
framework obtained by unfolding the polynomial
$$
f_{n,k}(z_1,\dots ,z_{n-1}):=\frac{1}{h}\,\sum_{i=1}^{n-1} y_i^h
$$
where $y$ and $z$ are related by
$$
\prod_{j=1}^{n-1}(1+y_jT)=1+\sum_{l=1}^{n-1}z_lT^l.
$$
In particular, if one assigns to  $z_l$ weight $l$,
$f_{n,k}$ becomes quasi--homogeneous of weight $h.$
Its unfolding space is spanned by the
classes of appropriate monomials, and a Zariski open
dense subset $V_{n,k}$ of this space carries the structure of
the Frobenius manifold as above. This subspace contains
the point $m$ corresponding to the {\it fusion potential}
$$
g_{n,k}(z_1,\dots ,z_{n-1}):=f_{n+1,k-1}(z_1,\dots ,z_{n-1},1).
$$
As D.~Gepner ([Ge2]) proved, the tangent space $\Cal{T}_mV_{n,k}$
with $\circ$-multiplication, that is, the
Milnor algebra of $g_{n,k},$ is isomorphic to the
Verlinde algebra (fusion ring) of the $su(n)_k$ WZW model
of the conformal field theory. Zuber in [Z] conjectured,
and Varchenko and Gusein--Zade in [G--ZV] proved, that the
lattice of the Verlinde algebra and the respective
bilinear form can be interpreted in terms of
vanishing cycles of $f_{n,k}.$

\smallskip

The total Frobenius manifold $V_{n,k}$ is thus
a deformation of this fusion ring, in much the same way
as quantum cohomology is the deformation of
the usual cohomology ring.

\medskip

{\bf 2.4. Direct sum diagram.} We will consider now
three Saito's frameworks $(p :\,N\to M, F, \omega )$
and $(p_i :\,N_i\to M_i, F_i, \omega_i ),\,i=1,2.$
We will call {\it the direct sum diagram}
any cartesian square
$$
\CD
N_1\times N_2 @>\nu>>     N \\
   @V(p_1,p_2)VV        @VVpV  \\
M_1\times M_2 @>\nu_M>> M    \\
\endCD
\eqno(2.8)
$$
with the following properties:

\smallskip

(i) $\nu^*(F)=F_1\boxplus F_2.$ 

\smallskip

(ii) $\nu^*(\omega )=\omega_1\boxtimes \omega_2.$

\medskip

Thus in a neighborhood of any point
of $N$ lying over the image of $\nu_M$ there exist
local coordinates $(z_a^{(1)},z_b^{(2)};t_e)$
such that $t_e$ are lifted from $M,$ and (i) can be written as
$$
F(z_a^{(1)},z_b^{(2)};\nu_M^*(t_e))=
F_1(z_a^{(1)};t_c^{(1)})+F_2(z_b^{(2)};t_d^{(2)})
\eqno(2.9)
$$
and similarly (ii) can be written as
$$
\omega(z_a^{(1)},z_b^{(2)};\nu_M^*(t_e))=
\omega_1(z_a^{(1)};t_c^{(1)})\wedge\omega_2(z_b^{(2)};t_d^{(2)}).
\eqno(2.10)
$$

\medskip

{\bf 2.4.1. Properties of the direct sum diagrams.}
Clearly, $\nu^{-1}(C)$ is defined by the equations
$d_{p_1}(F_1)\boxplus d_{p_2}(F_2)=0.$ Both summands
then must vanish so that $\nu^{-1}(C)=C_1\times C_2.$
Denote by $\nu_C:\,C_1\times C_2\to C$ the restriction of $\nu .$
From (2.9) one then sees that
$$
\nu_C^*(\roman{Hess}\,(F))=\roman{Hess}\,(F_1)
\boxtimes \roman{Hess}\,(F_2)
\eqno(2.11)
$$
and hence $\nu_M^{-1}(G)=G_1\times M_2\cup M_1\times G_2.$
Let now $m=\nu_M(m_1,m_2),\,m_i\in M_i.$ Choose  small neighborhoods
$m\in U$ in $M$, $m_i\in U_i$ in $M_i$ such that $\nu_M(U_1\times U_2)\subset U.$ Number the connected components $U_i^{(1)}$ of $p_{C_1}^{-1}(U_1),$
resp. $U_j^{(2)}$ of $p_{C_2}^{-1}(U_2),$ by some indices $i,$ resp. $j,$
as in (2.1). Then the connected components of $p_C^{-1}(U)$
are naturally numbered by the ordered pairs $I=(ij)$ in such a way
that 
$$
\nu_C (U_i^{(1)}\times U_j^{(2)})\subset U_I.
$$
From now on we will assume that all the frameworks we are considering
are strong ones. Then one can define $e_I,u^I,\eta_I$ etc as above,
and  from (2.3), (2.4), (2.9)--(2.11)
one immediately sees that
$$
u^I(m)=u_1^i(m_1)+u^j_2(m_2),\ \eta_I(m)=\eta^{(1)}_i(m_1)\eta^{(2)}_j(m_2)
\eqno(2.12)
$$
where in the right hand side we have the respective 
local functions on $M_1,M_2.$

\smallskip

The following slightly less evident restriction formula will be also needed
in the next section.

\medskip

\proclaim{\quad 2.4.2 Proposition} Let $I=(ij), K=(kl),
\eta_{IK}=e_I\eta_K=e_K\eta_I,$ and similarly
$\eta_{ik}^{(1)}=e_i\eta_k^{(1)}$ etc. Then we have
in the same notations as in (2.12):
$$
\eta_{IK}(m)=\delta_{jl}\eta_{ik}^{(1)}(m_1)\,\eta_l^{(2)}(m_2)+
\delta_{ik}\eta_{k}^{(1)}(m_1)\,\eta_{jl}^{(2)}(m_2).
\eqno(2.13)
$$
\endproclaim

\smallskip

{\bf Proof.} Calculate $\nu_M^*(d\eta_I)$ in two ways.
On the one hand, we have
$$
\nu_M^*(d\eta_I)=\sum_K \nu_M^*(\eta_{IK})\nu_M^*(du^K).
\eqno(2.14)
$$
As at the end of the proof of Prop. 2.2.2, we can
identify $du^K$ with a map from $p_{C*}(\Cal{O}_C)$ to
$\Cal{O}_M$ vanishing on all components except for
the $K$--th one where it is the canonical pushforward.
After restriction to $M_1\times M_2$ it may therefore
be non--vanishing only on $U_k^{(1)}\times U_l^{(2)}$
so that we can calculate $\nu_M^*(\eta_{IK})$ by restricting
$\nu_M^*(d\eta_I)$ to this product.

\smallskip

On the other hand, in view of (2.12),
$$
\nu_M^*(d\eta_I)=d\nu_M^*(\eta_I)=d(\eta_i^{(1)}\boxtimes
\eta_j^{(2)})=d\eta_i^{(1)}\boxtimes \eta_j^{(2)} +
\eta_i^{(1)}\boxtimes d\eta_j^{(2)}=
$$
$$
\sum_r\eta_{ir}^{(1)}du^r_1\boxtimes\eta_j^{(2)}+
\sum_s\eta_{i}^{(1)}\boxtimes\eta_{js}^{(2)}du^s_2.
\eqno(2.15)
$$
Only the $k$--th summand in the first sum
restricted to $U_k^{(1)}\times U_l^{(2)}$ may be non--vanishing
and considered as a map (cf. above) it
equals $\delta_{jl}\eta_{ik}^{(1)}\boxtimes\eta_j^{(2)}$
times the pushforward map. We have the similar
expression for the $l$--th summand of the second sum.
Comparison with (2.14) furnishes (2.13) because $\nu_M^*(du^K)=du_1^k\boxplus du_2^l.$

\medskip

{\bf 2.5. Direct sums of singularities.} In the theory of
singularities, we can compare the miniversal unfolding spaces $M_f,M_g,M_{f+g}$
of the germs $f$, $g$ and $f+g.$ It so happens
that they fit into the direct sum diagram
(2.8) (the only choice that remains is that
of the  volume form $\omega$ on the space of $f+g$
which is natural to take decomposable as in 2.4 above).

\smallskip

By iteration, we can consider arbitrary number of summands.
In particular, the Frobenius manifold $A_{n_1,\dots ,n_k}$
which is obtained by unfolding the
quasi--homogeneous singularity at zero $f(z):=z_1^{n_1}+\dots
+z_k^{n_k}$ is related to the summands $A_{n_i}$
in the way described above. We will show in the
next section, that the formulas (2.12) and (2.13)
imply a much neater description: $M_{f+g}=M_f\otimes M_g$,
and in particular $A_{n_1,\dots ,n_k}$
is the tensor product of $A_{n_i}$ in the sense of [KMK]
(in the context of formal
Frobenius manifolds) and [K] (in the global context).

\newpage

\centerline{\bf \S 3. Tensor products and their submanifolds}

\bigskip

{\bf 3.1. Tensor product of formal Frobenius manifolds.}
Let us first of all recall the general construction
of the tensor product of formal Frobenius manifolds
over a common coefficient ring $k$. Instead of
$(M=\roman{Spf}\, k[[H^t]],g ,\Phi )$ as
in the Definition 1.4.1 we will be writing 
$(H,g,\Phi ).$ We will not assume that our manifolds are
of qc--type.

\smallskip

Let $(H^{(i)},g^{(i)},\Phi^{(i)}),\,i=1,2,$ be two formal
Frobenius manifolds. Then
$$
(H^{(1)},g^{(1)},\Phi^{(1)})\otimes (H^{(2)},g^{(2)},\Phi^{(2)})=
(H^{(1)}\otimes H^{(2)},g^{(1)}\otimes g^{(2)},\Phi )
$$
where the terms of $\Phi$ of degree $n$ are defined in the following way: 
reconstruct
the Cohomological Field Theories $I^{(i)}_n$ for the two factors as in 1.5,
put $I_n:=I_n^{(1)}\cup I_n^{(2)}$ (cup product
in $H^*(\overline{M}_{0,n}))$ and cap $I_n$ with the fundamental
class of $\overline{M}_{0,n}.$

\smallskip

If the factors are endowed with flat identities
and Euler fields, with $d_0^{(i)}=1$ for both of them,
one can define in a canonical way the flat identity
and the Euler field for the product ([K]). Moreover, the
$d$--spectrum of the product is
$$
d=d^{(1)}+d^{(2)},\,\,\{q_A\}=\{q_a+q_b\}
$$
as a sum of families with multiplicities. Notice
that if $\roman{ad}\,E$ is semisimple for both factors,
it is semisimple for the product, and for $A=(a,b)$ one can identify 
$(\partial_A, q_A)$ with
$(\partial_a\otimes\partial_b,q_a+q_b).$

\smallskip

Since the tensor product potential is defined coefficientwise
and involves multiplication in all cohomology
algebras $H^*(\overline{M}_{0,n}),\,n\ge 3,$
it is practically impossible to calculate it
directly. The problem becomes much more manageable
if we deal with (germs of) semisimple analytic manifolds.

\smallskip

The reason for this is that generally a germ of
Frobenius manifold of a given dimension
depends on functional parameters, even in
the presence of a flat identity and Euler field.
The semisimplicity condition puts sufficiently
strong constraints in order to leave undetermined 
only a finite number
of constants, and then it is reasonable to rxpect
that the tensor product is calculable in terms of these
constants.
Below we review the relevant results
following [M] and [MM].

\medskip

{\bf 3.2. Moduli space and tensor product of germs of
semisimple Frobenius manifolds.} Consider a pointed germ $M$
of analytic Frobenius manifold over $\bold{C}$, (or a formal manifold
with zero as the base point), as usual, with flat $e$, Euler $E$
and $d_0=1,$ and having pure even dimension. It will be called 
{\it tame semisimple} if the operator $E\circ$ has simple
spectrum $(u^1_0,\dots ,u^n_0)$ on the tangent space to the base point.
We have the following general facts already partly invoked in the
specific situation of Saito's framework in 2.2 above:

\medskip

a) In a neighborhood of the base point, eigenvalues
$(u^1,\dots ,u^n)$ of $E\circ$ on $\Cal{T}_M$ form a
local coordinate system (Dubrovin's canonical
coordinates), taking the values $(u^1_0,\dots ,u^n_0)$
at the base point.
 The potential $\Phi$ is an analytic function of these
coordinates. If the initial manifold was only assumed 
to be formal, from tame semisimplicity it follows
that it is in fact the completion of a pointed analytic
germ.

\smallskip

b) Put $e_i=\partial /\partial u^i.$ Then $e_i\circ e_j=\delta_{ij}.$
In particular, $e=\sum_i e_i.$ It follows that the $\circ$ multiplication
on the tangent spaces is semisimple.

\smallskip

c) We have $g(e_i,e_j)=0$ for $i\ne j.$ Furthermore, there exists
a function $\eta$ defined up to addition of a constant
such that $g(e_i,e_i)=e_i\eta :=\eta_i.$ Moreover, we have
$eg =\roman{const}$, $Eg =(D-1)\eta +\roman{const}.$
Finally, $E=\sum_i u^ie_i.$

\medskip

A very important feature of canonical coordinates
is that a given tame semisimple germ can be uniquely
extended to the Frobenius structure on the universal
covering of the total $(u^i)$--space with deleted
partial diagonals. This follows from the Painlev\'e
property of the solutions of Schlesinger's equations:
cf. [M], Ch. II, sec. 1--3.
We will call this extension {\it the maximal tame continuation}
of the initial germ. The qualification
``tame'' is essential.
It may well happen that a further extension containing non--tame
semisimple points or even points with non--semisimple multiplication on
the tangent space is possible: e.g. points in $A_n$
where $F^{\prime}(z)$ has multiple roots have the latter
property.

\smallskip

\proclaim{\quad 3.2.1. Definition} Special coordinates
of a tame semisimple pointed germ of Frobenius manifold
consist of the values at the base point
of the following functions:
$$
(u^i,\,\eta_j, v_{ij}:=\frac{1}{2}\,(u^j-u^i)\,\frac{\eta_{ij}}{\eta_j}\,)
\eqno(3.1)
$$
Here $\eta_{ij}:=e_ie_j\eta$.
\endproclaim

\smallskip

To avoid any misunderstanding, let us stress that the canonical coordinates
are functions on a germ, whereas special coordinates
are functions on the moduli space of germs.

\smallskip

For a description of the necessary and generically sufficient
conditions for a system of mumbers to form
special coordinates of a Frobenius germ, see [MM], 2.7, pp. 26--27,
and 2.6, p. 23, where some inaccuracies of [M] are corrected.
The following Theorem summarizes the properties of
special coordinates that we will use. 
 
\medskip
\proclaim{\quad 3.2.2. Theorem} (i) Any tame semisimple pointed germ
with labelled spectrum of $E\circ$
is uniquely (up to isomorphism) defined by its special coordinates.

\smallskip

(ii) Let  $(u^{\prime i},\,\eta_j^{\prime}, v_{ij}^{\prime})$ for $i\in S$
and $(u^{\prime\prime i},\,\eta_j^{\prime\prime}, v_{ij}^{\prime\prime})$
for $j\in T$
be special coordinates of two pointed germs.
If the family $u^{\prime i}+u^{\prime\prime j}$ consists
of pairwise distinct elements, then the tensor product of the two germs defined
through their completions is again a tame semisimple pointed germ
whose canonical coordinates are naturally labelled by the
pairs $I\in S\times T$ and have the following form:
for $I=(i,j),\,K=(k,l),$
$$
u^I=u^{\prime i}+u^{\prime\prime j},\,
\eta_I=\eta_i^{\prime}\eta_j^{\prime\prime},\,
v_{IK}=\delta_{jl}v_{ik}^{\prime}+\delta_{ik}v_{jl}^{\prime\prime}.
\eqno(3.2)
$$
\smallskip

(iii) Let another two germs be obtained from the initial ones by
analytic continuation and subsequent shifts of base points.
Then their tensor product can be obtained from
the initial tensor product by analytic continuation and the
appropriate shift of the base point. In this sense,
the tensor product does not depend on the choice of
base points.
\endproclaim

\smallskip

The first statement is proved in [M] and [MM].
The second and the third ones are due to R.~Kaufmann ([K]).
Actually, the third statement is proved in [K] in
the considerably more general context: Kaufmann uses
flat coordinates and does not assume
semisimplicity or absence of odd coordinates.
The fact that the tensor product of two convergent
germs is again convergent is proved in his paper
in preparation,  without
semisimplicity assumption as well.

\smallskip

 Kaufmann remarks that in order to prove (3.2)
it suffices to control the relevant potentials
only to the fourth order in flat
coordinates, and the necessary calculation
can then be done directly.

\smallskip

We can now deduce from (3.2) the following corollary.

\medskip

\proclaim{\quad 3.2.3. Theorem} Assume that we have the
direct sum diagram of Saito's frameworks as in 2.4
above. Then the Frobenius manifold $M$ is
(canonically isomorphic to) the tensor product
of the Frobenius manifolds $M_1\otimes M_2.$
\endproclaim

\smallskip

{\bf Proof.} In the notation of (2.12) we may assume that
$m_1, m_2$ and $m$ are tame semisimple, because tameness
is the open property. Then (2.12) coincides with the
first two formulas of  (3.2). The third one
follows directly from  (2.13) and the definition of $v_{ij}$
in (3.1).

\medskip

We will now prove that the integral part of the spectrum
corresponds to a Frobenius submanifold.

\medskip

\proclaim{\quad 3.3. Proposition} Assume that we have
an analytic or formal Frobenius manifold $M$ with an Euler field
$E,\,d_0=1,\,D\in\bold{Z}$, and flat identity. Let $-\roman{ad}\, E$ be semisimple
on flat vector fields with spectrum $d_a$,
$(x_a)$ a flat coordinate system
with 
$$
E=\sum_{a:\,d_a\ne 0}d_ax^a\partial_a+\sum_{b:\,d_b=0}r_b\partial_b.
$$
and $e=\partial_0.$
Define the submanifold
$HM\subset M$ by the equations 
$$
x_c=0\ \roman{for\ all\ } c\ \roman{such\ that\ } d_c\notin \bold{Z}.
$$
Finally, assume that at least one of the following conditions
is satisfied:

\smallskip

(i) $r_b=0$ for all $b$ with $|\Delta_b|=2.$

\smallskip

(ii) $M$ is of qc--type, and $\sum_{b:\,d_b=0}r_b\Delta_b$
takes only integral values on $B$.

\smallskip

Then $HM$ with induced metric, $\circ$--multiplication, $E$
and $e$ is a Frobenius manifold.
\endproclaim

\smallskip

{\bf Remark.} From the proof it will be
clear that one can replace integers in this statement
by any arithmetic progression containing 0
to which $D$ and $d_0$ belong.

\medskip

{\bf  Proof of the proposition 3.3.} If $d_c$ is not integral,
the functions
$Ex_c=d_cx_c,\,ex_c=0$ vanish on $HM$. Hence $E$ and $e$
are tangent to $HM$ and can be restricted to it. From the equation
$(d_a+d_b-D)g_{ab}=0$ ([M], p. 32, (2.17)) one sees that
if $d_a,\,D\in \bold{Z},\,d_b\notin \bold{Z},$ we have $g_{ab}=0.$
Therefore the restriction of $g$ to $HM$ is non--degenerate
(it is obviously flat), and $x_a$ for $d_a\in \bold{Z}$ restrict to a flat
coordinate system on $HM.$ The $\circ$--product of
two vector fields tangent to $HM$ at the points of $M$
does not contain the transverse components. In fact,
we have $E\Phi_{ab}{}^c=(d_0-d_a-d_b+d_c)\Phi_{ab}{}^c$
([M], p. 32, (2.18)).
Hence if $d_a,d_b\in\bold{Z}$, $d_c\notin\bold{Z}$,
then in the case (i)
every monomial in the series $\Phi_{ab}^c$ must
be an eigenvector of $E$ with non--integral eigenvalue,
and therefore it must 
contain some $x_e$ with $d_e\notin\bold{Z}$ so that
it vanishes on $HM.$ In the case (ii) we apply the same reasoning
separately to the generalized (involving exponentials)
 monomials contributing to the third
derivatives of $\Psi$ in (1.9)
and to the third derivatives of $c$.

\smallskip

The same reasoning shows that the induced multiplication of
vector fields on $\Cal{T}_{HM}$ is defined by the third
derivatives of the induced potential.

\medskip

{\bf 3.4. Special coordinates of $A_n$.} We return to the
notation of 2.3.1.

\medskip

\proclaim{\quad 3.4.1. Proposition} Consider the points of $A_n$ where
$a_1,\dots ,a_{n-2}=0,\, a_{n-1},\, a_n$ arbitrary.
Choose a primitive root $\zeta$ of $\zeta^n=1$ and
a root $b$ of $b^n=-\dfrac{a_{n-1}}{n+1}.$
At these points we have:
$$
u^i=a_n+\frac{n}{n+1}\zeta^ia_{n-1}b,
\eqno(3.3)
$$
$$
\eta_i=\frac{\zeta^i}{n(n+1)\,b^{n-1}},
\eqno(3.4)
$$
$$
v_{jk}=\frac{1}{(n+1)\,(1-\zeta^{k-j})}.
\eqno(3.5)
$$
\endproclaim

\smallskip

 {\bf Remark.} It is
suggestive to compare these coordinates
with those for the quantum cohomology
of $\bold{P}^{n-1}$ ([M], p. 71) on the plane
spanned by the identity (coordinate $x_0$) and the dual
hyperplane section (coordinate $x_1$):
$$
u^i=x_0+n\,\zeta^ie^{\frac{x_1}{n}},
\eqno(3.3a)
$$
$$
\eta_i=\frac{\zeta^i}{n}\,e^{-x_1\,\frac{n-1}{n}},
\eqno(3.4a)
$$
$$
v_{jk}=\frac{1}{1-\zeta^{k-j}}.
\eqno(3.5a)
$$

\medskip

{\bf Proof.} At our subspace $F(z)=z^{n+1}+a_{n-1}z+a_n.$
Hence 
$$
F^{\prime}(z)=(n+1)(z^n+\dfrac{a_{n-1}}{n+1})
$$
has roots $\rho_i=\zeta^ib.$
But for $A_n$--manifolds we have universally
$u^i=F(\rho_i),\,\eta_i=\dfrac{1}{F^{\prime\prime}(\rho_i)}.$
 This furnishes (3.3) and (3.4).

\smallskip

The proof of (3.5) is longer. We have to calculate
the values of functions
$$
\frac{1}{2}\,(u^k-u^j)\,\frac{\eta_{jk}}{\eta_k},
\eqno(3.6)
$$
restricted to the plane of our base points.

\smallskip

At a generic point of $A_n$, we can calculate
$\eta_{jk}$ in the following three--step way:
$$
\eta_{jk}=\frac{\partial \eta_j}{\partial u^k}=
\sum_{l,m=1}^n \frac{\partial \eta_j}{\partial \rho_m}\,
\frac{\partial \rho_m}{\partial a_l}\,
\frac{\partial a_l}{\partial u^k}.
\eqno(3.7)
$$
Since
$$
\eta_j=\frac{1}{F^{\prime\prime}(z)}=
\frac{1}{(n+1)\prod_{i:\,i\ne j}(\rho_i-\rho_j)},
$$
we have
$$
\frac{\partial \eta_j}{\partial \rho_m}=-\frac{\eta_j}{\rho_m-\rho_j},
\eqno(3.8)
$$
if $m\ne j$, and
$$
\frac{\partial \eta_j}{\partial \rho_j}=
\eta_j\sum_{i:\,i\ne j}\frac{1}{\rho_i-\rho_j}.
\eqno(3.9)
$$
Moreover,
$$
\frac{\partial \rho_m}{\partial a_l}=-(n-l)\,\rho_m^{n-l-1}\eta_m.
\eqno(3.10)
$$
This can be checked by derivating the identity $F^{\prime}(\rho_m)=0.$

\smallskip

Finally, according to [M], p. 47, (4.24), we have
$$
\sum_{l=1}^n\frac{\partial a_l}{\partial u^k}\,\rho_i^{n-l}=\delta_{ik}.
\eqno(3.11)
$$
We will now restrict (3.8)--(3.11) to our plane.

\smallskip

Using (3.3) and (3.4), we  get consecutively:
$$
\frac{\partial \eta_j}{\partial \rho_m}=\frac{1}{n\,a_{n-1}}\,
\frac{1}{\zeta^{m-j}-1}
\eqno(3.12)
$$
if $m\ne j$, and
$$
\frac{\partial \eta_j}{\partial \rho_j}=\frac{n-1}{2n\,a_{n-1}},
\eqno(3.13)
$$
$$
\frac{\partial \rho_m}{\partial a_l}=
-\frac{n-l}{n(n+1)}\,b^{-l}\zeta^{-ml}.
\eqno(3.14)
$$
Solving (3.11) for partial derivatives, we also find
$$
\frac{\partial a_l}{\partial u^k}=\frac{1}{n}\,b^{l-n}\zeta^{kl}.
\eqno(3.15)
$$
It remains to substitute (3.12)--(3.14) into (2.7) to get after some
calculation
$$
\eta_{jk}=\frac{2\zeta^{k-j}}{(\zeta^{k-j}-1)^2}\,
\frac{1}{n\,a_{n-1}^2}.
\eqno(3.16)
$$
Finally, substituting (3.3), (3.4) and (3.16) in (3.6),
we obtain (3.5).

\medskip

{\bf 3.5. Tensor products of $A_n$'s.} We want to describe
$(n_1,\dots ,n_N)$ with non--trivial $H(A_{n_1}\otimes \dots\otimes
A_{n_N}).$ We can assume $n_i\ge 2$ because $A_1$
is identity with respect to the tensor multiplication.
The first necessary condition, following from
(1.20) is
$$
\bold{d}:=\sum_{i=1}^N \frac{n_i-1}{n_i+1}\in \bold{Z}.
\eqno(3.17)
$$
If it is satisfied, the full $d$--spectrum of the tensor
product consists of certain rational points between $0$ and
$\bold{d}.$ Multiplicity of $0$ and $\bold{d}$ is one.
Generally, the multiplicity of some $m\le \bold{d}$
is
$$
h^{2m}(H(\otimes_i A_{n_i})):=
\roman{the\ number\ of\ } (i_1,\dots ,i_N)\in\bold{Z}^N_{\ge 0}\
\roman{satisfying}
$$
$$
\sum_{k=1}^N\frac{i_k}{n_k+1}=m,\, 0\le i_k\le n_k-1.
\eqno(3.18)
$$
The $d$--spectrum of $H(\otimes A_{n_k})$ consists of the
part of (3.18) for  all integer $m$. Clearly, (0.1)
is a particular case of (3.18).

\smallskip

The flat part of $E$ in the total tensor product
and in the $H$--part of it vanishes because it
vanishes on all factors (cf. [K], Theorem 6.3).

\smallskip

Let us show that $\otimes_kA_{n_k}$ admits a tame semisimple
base point which is the ``sum'' of the points 
the special coordinates of which we have calculated.
Choose base points on all $A_{n_k}$ as in Prop. 3.4.1. 
For our purpose, we may even assume that $a_{n_k}=0$
on each $A_{n_k}.$ Therefore, slightly changing notation
of (3.3), we will assume that canonical coordinates
of the base point of $A_{n_k}$ are of the form $u^i=\zeta_k^ic_k,$
where $\zeta_k$ is the primitive root of unit of degree $n_k$
and $0\le i\le n_k-1.$ Then in view of the Theorem 3.2.2 for
$I=(i_1,\dots ,i_N),\, J=(j_1,\dots ,j_N)$ we have
a pair of canonical coordinates $u^I,u^J$ on the full
tensor product whose values at the base point
of this product are
$$
u^I=\sum_{k=1}^N\zeta_k^{i_k}c_k,\,\,u^J=\sum_{k=1}^N\zeta_k^{j_k}c_k.
$$
One easily sees that with generic choice of $c_k$ these
coordinates are distinct for all $I\ne J.$

\medskip

{\bf 3.6. Involutive pairs of Gepner's manifolds.}
In the notations of 2.3.2, consider a pair of Frobenius manifolds
$V_{n+1,k}$ and $V_{k+1,n}.$ They contain respectively
points producing Verlinde's algebras of $su(n+1)_k$
and $su(k+1)_n$ which are isomorphic according to the
reasoning of [G--ZV] (the level--rank duality). This reasoning runs as follows.
Assume for concreteness that $n\ge k+1.$ Consider the function
$$
\widetilde{f}_{n+1,k}(z_1,\dots ,z_n):=f_{k+1,n}(z_1,\dots ,z_k)+
\sum_{i=k+1}^nz_iz_{n+k-i+1}.
\eqno(3.19)
$$
Then $f_{n+1,k}$ and $\widetilde{f}_{n+1,k}$ are 
quasi--homogeneous polynomials of the same degree,
depending on the same set of weighted variables, and
having an isolated critical point at the origin.
Hence they belong to a connected family of polynomials
with the same property, and whatever structures
can be derived from their lattices of vanishing cycles,
they can be identified. On the other hand,
$\widetilde{f}_{n+1,k}$ is obtained from
$f_{k+1,n}(z_1,\dots ,z_k)$ by adding a sum of squares
which again does not change the structure of vanishing cycles,
except for that of the intersection form
which changes in a controlled way.

\smallskip

In fact, adding a sum of squares does not
change the respective Frobenius manifolds:
this agrees with the fact that $A_1$ is the tensor identity.
Hence $V_{n+1,k},$ or at least its germ at the origin,
is deformable to (the germ of) $V_{k+1,n}.$ 

\smallskip

Perhaps, these Frobenius manifolds, 
or at least their appropriate coverings, are themselves
isomorphic. To check this, it would suffice
to identify their special coordinates at an appropriate
pair of tame semisimple points. The Gepner--Verlinde
points $m\in V_{n,k}$ (cf. 2.3.2 above) are certainly not tame
because  the potential has only $n$, not $n+k$
different critical values at $m$.

\newpage

\centerline{\bf \S 4. Maurer--Cartan and master equations}

\bigskip

{\bf 4.1. Maurer--Cartan equations.} 
Fix a supercommutative $\bold{Q}$--algebra $k$.
All our
structures are $\bold{Z}_2$--graded, notation like
$\widetilde{x}$ means the parity of a homogeneous element
$x.$ Let $g=g_0\oplus g_1$ be a Lie superalgebra over $k$,
supplied with an odd differential $d$ satisfying
$d[a,b]=[da,b]+(-1)^{\widetilde{a}}[a,db].$ 

\smallskip

Put $Z=Z(g,d):=\roman{Ker}\,d, B=B(g,d):=\roman{Im}\,d, H:=H(g,d)=Z/B.$
Clearly, $Z$ is a Lie subalgebra, and $B$ its ideal, so that
$H$ with induced bracket product is a Lie superalgebra.

\smallskip

The differential $d$ can be shifted.

For $\gamma \in g_1,$
put $d_{\gamma}(a):=da +[\gamma ,a].$ Clearly,
$d_{\gamma}[a,b]=[d_{\gamma}a,b]+(-1)^{\widetilde{a}}[a,d_{\gamma}b].$

\medskip

\proclaim{\quad 4.1.1. Claim} a) We have $d_{\gamma }^2=0$ if
$$
d\gamma +\frac{1}{2}[\gamma ,\gamma]=0.
 \eqno(4.1)
$$

b) Let $\gamma^{\prime}=\gamma +\epsilon\beta$, $\beta\in g,$
$\epsilon$ an even or odd constant with $\epsilon^2 =0$ 
such that $\epsilon\beta$ is odd. Assume that $\gamma$  
satisfies (4.1). Then  $\gamma^{\prime}$ satisfies (4.1) as well iff
$$
d_{\gamma}(\beta )=d\beta +[\gamma ,\beta ]=0.
 \eqno(4.2)
$$
\endproclaim

\smallskip

This is straightforward.

\smallskip

If $K$ is another supercommutative $k$--algebra,
we  define  $g_K=K\otimes_kg, d_K=1\otimes d.$
We will always work with $K$ flat over $k$ so that
$Z_K:=\roman{Ker}\,d_K=K\otimes_kZ$, and similarly for 
$B$ and $H.$ Claim 4.4.1 is of course applicable to
$(g_K,d_K)$ as well.

\medskip

We want to produce from this setting a
non--linear version of the homology $H(g,d)$
or rather of the diagram $g\supset Z\to H.$

\smallskip

The most straightforward is the case when,
say, $g$ is free of finite rank over $k.$
We then replace $g$ by the linear superspace
$\Cal{G}:=\roman{Spec}\,k[\Pi g^t],$
where $\Pi$ is the parity inversion functor,
$Z$ by the closed subspace $\Cal{Z}\subset \Cal{G}$
defined by the equations (4.1). In order to understand what
should be the non--linear version of $B$, we
interpret the Claim 4.1.1 (ii) as saying that
$d_{\gamma}$--cycles form the Zariski tangent space
to the point $\gamma$ of the Maurer--Cartan space (4.1).
It then contains the subspace of $d_{\gamma}$--boundaries,
and we can construct the distribution $\Cal{B}$
generated by the boundaries.
If the quotient space $\Cal{H}=\Cal{Z} /\Cal{B}$
in some sense exists, it can be regarded as 
the non--linear cohomology of $(g,d)$.

\smallskip

In more down--to--earth terms,
choose a (homogeneous) basis $\{\gamma_i\}$ of $g$ and a family
of independent (super)commuting variables $t^i$ such that
$\widetilde{t}^i=\widetilde{\gamma}_i+1.$
Then $\Gamma :=\sum_it^i\gamma_i$ is a generic odd element
of $g$ (or rather of $k[t^i]\otimes g$), and the equation
$d\Gamma +\dfrac{1}{2}[\Gamma ,\Gamma ]=0$ is equivalent to the
system of equations
$$
\forall k:\quad \sum_i(-1)^it^iD_i^k+\frac{1}{2}\sum_{i,j}
t^it^j(-1)^{i(j+1)}L_{ij}{}^k=0.
 \eqno(4.3)
$$
Here we define the structure constants by
$d\gamma_i=\sum_kD_i^k\gamma_k$ and $[\gamma_i,\gamma_j]=
\sum_kL_{ij}{}^k\gamma_k$ and use the following
shorthand for the signs: $(-1)^{i(j+1)}$ means
$(-1)^{\widetilde{\gamma}_i(\widetilde{\gamma}_j+1)}$ etc.

\smallskip

These equations define the coordinate ring $R$ of the
affine scheme which we called $\Cal{Z}.$ Obviously, $\Cal{Z}$ represents the following
functor on the category of supecommutative $k$--algebras $K$:
$$
K\mapsto \{solutions\ to\ (4.1)\ in\ (K\otimes g)_1\}.
 \eqno(4.4)
$$ 

\smallskip

Similarly, if we have any odd $d_{\Gamma}$--cycle
$\epsilon\beta =\sum_a\epsilon s^a\gamma_a$ with coefficients 
in $K\otimes R,$ 
the statement 4.1.1 b) means that the map
$X_{\beta}:\,t^a\to s^a$ descends to the derivation
of $K\otimes R$ over $K$ that is, to a vector field
on $\Cal{G}_K$ of parity $\widetilde{\epsilon}.$
Of course, the adequate functorial language for
derivations is that of the first order 
infinitesimal deformations of points, because generally
the vector fields implied by 4.1.1 b) are defined
only in the infinitesimal neighborhood of $\gamma$.

\smallskip

We will stop now discussing the case of finite
rank $g$ because in most interesting examples
this does not hold, and only $H(g,d)$
is of finite rank. 

\smallskip

So we step back and try to produce
a formal section of $\Cal{Z}$ passing through $\gamma =0$ 
and transversal
to the distribution $\Cal{B}$. We want it
to be of the same size as $H$, or rather $\Pi H$,
and we will
assume henceforth that $H$ is free of finite rank.
From now on in this section, we denote  $K:= k[[\Pi H^t]]=k[[x_i]]$
where $x_i$ are coordinate functions
on $\Pi H$ dual to a basis of $\Pi H.$
Any element $\Gamma \in g_K$ can be uniquely written
as $\sum_{n\ge 0} \Gamma_n$ where $\Gamma_n$
is homogeneous of degree $n$ in $x_i$.
Such an element can be naturally called
a formal section of $\Cal{Z}$, or {\it a generic
(formal) solution to (4.1)}, if it has the following
properties:

\medskip

a) $\Gamma \in (g_K)_1, \Gamma_0 =0, \Gamma_1 = \sum_i x_ic_i$
where $dc_i=0$ and classes of $c_i$ form a basis
of $H$ odd dual to $\{x_i\}.$

\smallskip

b) $d_K\Gamma +\dfrac{1}{2}[\Gamma ,\Gamma ]=0.$ 

\smallskip

The necessary condition for the existence of $\Gamma$
is the identical vanishing of the Lie bracket induced on $H(g,d).$
In fact, the equation $d\Gamma +\dfrac{1}{2}[\Gamma ,\Gamma ]=0$
implies (assuming a) above) $d\Gamma_2 +\dfrac{1}{2}[\Gamma_1 ,\Gamma_1 ]=0.$
Hence $[c_i,c_j]\in B.$ However, generally it is
not sufficient. In fact, the next equation reads
$d\Gamma_3+[\Gamma_1,\Gamma_2]=0,$ but since $\Gamma_2$
may be non--closed, we cannot conclude that
$[\Gamma_1,\Gamma_2]$ is a boundary.
The manageable sufficient condition is stronger:
$(g,d)$ must be quasi--isomorphic to the differential
Lie algebra $H(g,d)$ with zero bracket and zero differential.
For a considerably more general treatment see [GoM]. 
Our direct and elementary approach is self--contained and
produces slightly more detailed information in the cases
essential for the theory of Frobenius manifolds.

\medskip

\proclaim{\quad 4.2. Theorem} (i) Assume that
there exists a surjective morphism of
differential Lie superalgebras $\phi: (g,[,],d)\to (H,0,0)$
inducing isomorphism on the homology.
Then there exists a generic formal solution $\Gamma$
to (4.1).

\smallskip

Moreover, $\Gamma$ can be chosen in such a way
that for any $n\ge 2,$ $\Gamma_n\in K\otimes \roman{Ker}\,\phi .$
In other words, $(\roman{id}\otimes \phi ) (\Gamma )=
\sum_ix_i[c_i].$ Such a solution will be called normalized.

\smallskip

(ii) If (i) is satisfied, then for any generic solution
$\Gamma ,$ non necessarily normalized, 
the map $\phi_K=\roman{id}\otimes\phi :\,
g_K\to H_K$ is the surjective morphism of
differential Lie superalgebras $ (g_K,[,]_K,d_{K,\Gamma})\to (H_K,0,0)$
inducing isomorphism on the homology.
\endproclaim

\medskip

{\bf Proof.} (i) Let $n\ge 1.$
Assuming that $\Gamma_i$ for $i\le n$ are already 
constructed, 
and writing $d$ instead of $d_K$ 
we must find $\Gamma_{n+1}$ from the equation
$$
d\Gamma_{n+1}=-\frac{1}{2}\sum_{i,j:\,i+j=n+1}[\Gamma_{i},\Gamma_j].
 \eqno(4.5)
$$
First of all we check that the right hand side of (4.5)
is closed in $g_K$. In fact, since the components $\Gamma_1,\ldots ,\Gamma_n$
satisfy the similar equations by the inductive assumption,
the differential of the rhs equals
$$
\frac{1}{2}\sum_{i+j+k=n+1}[[\Gamma_i,\Gamma_j],\Gamma_k].
$$
This expression vanishes because the Jacobi identity for
odd elements reads
$$
[[\Gamma_i,\Gamma_j],\Gamma_k]+[[\Gamma_k,\Gamma_i],\Gamma_j]
+[[\Gamma_j,\Gamma_k],\Gamma_i]=0.
$$
Hence the coefficients of the rhs of (4.5)
(as polynomials in $x_i$) belong to $Z\,\cap\, [g,g].$
But $[g,g]\in \roman{Ker}\,\phi$ and
$Z\cap \roman{Ker}\,\phi =B$ because $\phi$
is a quasi--isomorphism. Thus we can solve (4.5).

\smallskip

We can add to any solution elements of $Z_K$ of degree $n+1.$
But $Z+\roman{Ker}\,\phi =g$ because $\phi$ induces surjection
on homology. Hence we can normalize $\Gamma_{n+1}$ by
the requirement $\Gamma_{n+1}\in K\otimes \roman{Ker}\,\phi .$

\medskip

(ii) Now fix $\Gamma$ satisfying (4.5) for all $n.$
We will write $d_{\Gamma}$ instead of $d_{K,\Gamma}$
and put $Z_{\Gamma}:=\roman{Ker}\,d_{\Gamma}\subset g_K$,
$B_{\Gamma}=d_{\Gamma}(g_K).$ We have $B\subset \roman{Ker}\,\phi_K$
and $[g_K,g_K]\subset \roman{Ker}\,\phi_K$, hence
$B_{\Gamma}\subset \roman{Ker}\,\phi_K.$ Therefore, $\phi_K$
is compatible with zero bracket and zero differential on $H_K.$
The natural inclusion $Z_{\Gamma}+\roman{Ker}\,\phi_K \to g_K$
becomes surjection after the reduction modulo the
 ideal $(x_i)$ of $K$, because $\phi$ is  surjective.
Hence this inclusion is surjective, and $\phi_K$ is surjective as well. It remains to show
that $\phi_K$ induces injection on homology, that is,
$$
Z_{\Gamma}\cap\, \roman{Ker}\,\phi_K \subset B_{\Gamma}.
 \eqno(4.6)
$$
Let $c=\sum_{n\ge 0}c_n\in Z_{\Gamma}.$ This means that
$dc_0=0$ and in general
$$
dc_n=-\sum_{i+j=n}[\Gamma_i,c_j].
 \eqno(4.7)
$$
(we keep writing $d$ for $d_K$).
Assuming that $\phi_K(c)=0$ we want to deduce
the existence of homogeneous elements $a_n$ of degree
$n$ in $g_K$ such that
$$
c_{n+1}=da_{n+1}+\sum_{i+j=n+1}[\Gamma_i,a_j].
 \eqno(4.8)
$$
We have $dc_0=0$ and $\phi (c_0)=0$, hence
$c_0$ is a boundary because $\phi$ is the
quasi--isomorphism. Assuming that
$a_0,\ldots ,a_n$ are found, we will establish
the existence of $a_{n+1}$ satisfying (4.8),
if we manage to prove that
$c_{n+1}-\sum_{i+j=n+1}[\Gamma_i,a_j]$ is 
$d$--closed. In fact, this element also
belongs to $\roman{Ker}\,\phi_K$ and so must
be a boundary.

\smallskip

The differential of this element is
$$
dc_{n+1}+\sum_{i+j=n+1}[\Gamma_i,da_j] 
-\sum_{i+j=n+1}[d\Gamma_i,a_j].
\eqno(4.9)
$$
Replace in the first sum $da_j$ by
$c_j-\sum_{k+l=j}[\Gamma_k,a_l]$ for $j\le n$
(this holds by induction). Replace in the second
sum $d\Gamma_i$ by the sum of commutators from (4.5).
The terms containing $c_j$ will cancel thanks to
(4.7). The remaining terms can
be written as
$$
-\sum_{i+j+k=n+1}[\Gamma_i,[\Gamma_j,a_k]]+
\frac{1}{2}\sum_{i+j+k=n+1}[[\Gamma_i,\Gamma_j],a_k].
$$
This expression vanishes because of Jacobi identity.

\medskip

\proclaim{\quad 4.2.1. Corollary} Define the map
$\psi =\psi_{\Gamma}: H_K\to g_K$ as the $K$--linear extension of
$$
H\to g_K:\ X\mapsto \overline{X}\Gamma
$$
where $\overline{X}$ acts on $K\otimes g$ as the 
 right $g$--linear extension
of the derivation on $K$ acting as $(\Pi X,\,*)$ on $\Pi H^t.$

\smallskip

Then $\psi$ is a section of $\phi_K$ if $\Gamma$ is normalized.
\endproclaim

\smallskip

{\bf Proof.} First of all, we have 
$\overline{X}(d\Gamma +\dfrac{1}{2}[\Gamma ,\Gamma ])=0$
from which it follows that $d_{\Gamma} (\overline{X}\Gamma )=0,$
that is, $\overline{X}\Gamma$ is a $d_{\Gamma}$--cycle.
Its image in $H_K$ is $(\overline{X}\Gamma_1 + \sum_{n\ge 2}
\overline{X}\Gamma_n)\ \roman{mod}\,B_{\Gamma}.$
The first term is clearly $X.$ The remaining ones
are in $K\otimes \roman{Ker}\,\phi_K,$ if $\Gamma$
is normalized.

\medskip

{\bf 4.3. Odd Lie (super)algebras.} As in 4.1, let
now $g=g_0\oplus g_1$ be a $k$--module endowed
with a bilinear operation {\it odd bracket} $(a,b)\mapsto [a\b b]$
which satisfies the following conditions:

\smallskip

a) {\it parity of $[a\b b]$ equals $\widetilde{a}+\widetilde{b}
+1,$}

\smallskip

b) {\it odd anticommutativity:}
$$
[a\b b]=-(-1)^{(\widetilde{a}+1)(\widetilde{b}+1)}[b\b a],
 \eqno(4.10)
$$
\smallskip

c) {\it odd Jacobi identity:}
$$
\quad [a\b [b\b c]]=[[a\b b]\b c] +(-1)^{(\widetilde{a}+1)(\widetilde{b}+1)}[b\b[a\b c]].
 \eqno(4.11)
$$
Such a structure will be called an odd Lie (super)algebra.
We consider such algebras endowed with an odd differential
satisfying
$$
d[a\b b]=[da\b b]+(-1)^{\widetilde{a}+1}[a\b db].
 \eqno(4.12)
$$
Physicists sometimes denote such multiplication
$\{,\}$ (see e.g. [LZ]). Our choice of notation
allows one to use consistently the standard sign mnemonics of
superalgebra, if $\b$ counts as an element
of parity one. 

\smallskip

If $(g,d)$ is the usual differential Lie superalgebra,
the parity change functor $g\mapsto \Pi g$ 
turns the usual bracket product $[\,,]$ into the odd bracket product,
and defines an equivalence of the two categories
(the differential changes sign).
It seems therefore that there is not much point
in considering odd brackets. However, in the
context of GBV--algebras they come together with
usual supercommutative multiplication, and
parity change then turns this multiplication into
odd one (see the next section). This is, of course, a
particular case of the general operadic formalism
over the category of superspaces, where
any operation can be inherently even or odd.

\smallskip

In the next section we choose to work with even
multiplication and odd bracket product. 
But we will use the results of this section,
with appropriately modified parities and signs, for
odd Lie superalgebras. In particular, the odd
Maurer--Cartan equation in  
the physical literature is called {\it the master equation}:
$$
d\Gamma +\frac{1}{2}[\Gamma\b \Gamma]=0.
 \eqno(4.13)
$$
The Theorem 4.2 provides conditions of its
solvability in 
$k[[H^t]]\otimes g$ rather than $k[[\Pi H^t]]\otimes g$.
Notice also that $\Gamma$ in (4.13) must be even.

\newpage

\centerline{\bf \S 5. Gerstenhaber--Batalin--Vilkovisky algebras}

\bigskip

{\bf 5.1. Gerstenhaber--Batalin--Vilkovisky algebras.} Let $\A$ be a
supercommutative algebra with identity over another
supercommutative algebra $k$ (constants).
Consider an odd
$k$--linear operator $\Delta : \A \to \A , \D (1)=0,$ with the following
property:
$$
\forall a\in \A,\qquad \p_a :=(-1)^{\widetilde{a}}([\Delta ,l_a]-l_{\Delta a})
$$
$$
\roman{is\ the\ } \roman{derivation\ of\ parity\ } \widetilde{a}+1\ 
\roman{over}\ k.
 \eqno(5.1)
$$
Here $l_a$ denotes the operator of left multiplication by $a$,
and brackets denote the supercommutator.
Explicitly,
$$
\p_a{b}=(-1)^{\widetilde{a}}\D (ab)-(-1)^{\widetilde{a}}(\D a)\,b-
a\,\D b.
$$
The sign ensures the identity $\p_{ca}=c\p_a$ for any constant
$c$. By definition of derivation,
$$
[\p_a ,\l_b]=l_{\p_ab}.
 \eqno(5.2)
$$
The pair $(\A,\D )$ is called a {\it GBV--algebra} if, in addition,
$\D^2=0.$ There is an obvious operation of scalar extension.

\medskip

\proclaim{\quad 5.1.1. Lemma} In any GBV--algebra we have
$$
[\D ,\p_a]=\p_{\D a},
 \eqno(5.3)
$$
$$
[\p_a, \p_b]=\p_{\p_ab}
 \eqno(5.4)
$$
\endproclaim

{\bf Proof.} From (5.1) we have
$$
[\D ,\p_a]=(-1)^{\widetilde{a}}([\D,[\D,\l_a]]-[\D,l_{\D a}]).
$$
From the Jacobi identity for operators and $[\D ,\D ]=0$
we find $[\D,[\D,\l_a]]=0$ because
$$
[\D,[\D,\l_a]]=[[\D,\D],\l_a]-[\D,[\D,\l_a]])
$$
 From (5.1) with $\D a$ replacing $a$
we have $[\D ,l_{\D a}]=(-1)^{\widetilde{a}+1}\p_{\D a}.$ Hence 
$$
[\D ,\p_a]=(-1)^{\widetilde{a} +1}[\D,l_{\D a}]=\p_{\D a}.
$$
To prove (5.4), we notice that since $[\p_a, \p_b]$ must be
a derivation, in the intermediate calculations we are allowed not to register
all the summands which are left multiplications: they
will cancel anyway. So we have, denoting such summands by dots
and using consecutively (5.1), (5.2), Jacobi and (5.3), and again (5.1) with 
$\p_{a}b$ replacing $a$:
$$
[\p_a, \p_b]=(-1)^{\widetilde{b}}[\p_a, [\D,l_b]-l_{\D b}]=
(-1)^{\widetilde{b}}[\p_a, [\D,l_b]]+\ldots =
$$
$$
(-1)^{\widetilde{b}+\widetilde{a}+1}[\D,[\p_a,l_b]]+\ldots =
(-1)^{\widetilde{b}+\widetilde{a}+1}[\D,l_{\p_ab}]+\ldots =
\p_{\p_ab}.
$$

\medskip

Define now the odd bracket operation on $\A$ by the formula
$$
[a\b b]:=\p_ab.
 \eqno(5.5)
$$
\medskip

\proclaim{\quad 5.1.2. Proposition} The pair of bilinear operations
(multiplication and odd bracket) defines on $\A$ the structure
of the odd Poisson algebra in the following sense:

\smallskip

(i) The odd bracket satisfies the odd anticommutativity, the odd
Jacobi and the odd Poisson identities:
$$
[a\b b]=-(-1)^{(\widetilde{a}+1)(\widetilde{b}+1)}[b\b a],
\quad [a\b [b\b c]]=[[a\b b]\b c] +(-1)^{(\widetilde{a}+1)(\widetilde{b}+1)}[b\b[a\b c]],
$$
$$
[a\b bc]=[a\b b]c+(-1)^{\widetilde{a}(\widetilde{b}+1)}b\,[a\b c].
 \eqno(5.6)
$$

\smallskip

(ii) $\D$ is the derivation with respect to the odd brackets
so that $(\A,\b ,\D )$ is the differential odd Lie algebra.
 
\endproclaim

{\bf Proof.} The anticommutativity can be checked directly. 
The Jacobi identity follows from (5.4) written as 
$[\p_a, \p_b]=\p_{[a\b b]}.$  The Poisson identity means
that $\p_a$ is a derivation. 
The last statement follows from (5.1).

\smallskip

Notice that with respect to the usual multiplication
$\D$ is the differential operator of order $\le 2$ and
not necessarily derivation.

\bigskip

{\bf 5.2. Additional differential.} Assume now that
we have an additional $k$--linear odd map $\d :\A\to \A$
which is the derivation with respect to the multiplicative
structure of $\A$ satisfying
$$
\d^2=[\d ,\D]=\d\D +\D\d=0.
 \eqno(5.7)
$$
We will say that $(\A,\D ,\d )$ is {\it a differential GBV--algebra
(dGBV)}.
\medskip

\proclaim{\quad 5.2.1. Lemma} We have
$$
[\,\d ,\p_a]=\p_{\d a}.
 \eqno(5.8)
$$
Therefore $\d$ is the derivation with respect to the odd bracket
as well. 
\endproclaim

\smallskip

{\bf Proof.} Since $[\,\d ,\p_a]$ is a derivation of $\A$,
we can calculate omitting the multiplication operators as above:
$$
[\,\d ,\p_a]=(-1)^{\widetilde{a}}[\d,[\D ,l_a]]+\ldots =
(-1)^{\widetilde{a}}([[\,\d ,\D],l_a]-[\D,[\,\d ,l_a]])+\ldots =
$$
$$
-(-1)^{\widetilde{a}}[\D,l_{\delta a}]+ \ldots =\p_{\delta a}.
$$
Furthermore,
$$
\delta\,[a\b b]=\delta\,\p_ab=[\d ,\p_a]b+(-1)^{\widetilde{a}+1}\p_a\d b=
$$
$$
=\p_{\d a}b+(-1)^{\widetilde{a}+1}\p_a\d b=[\,\d a\b b]+
(-1)^{\widetilde{a}+1}[a\b\d b].
 \eqno(5.9)
$$

\medskip

{\bf 5.2.2. Shifted differential.} Let $\d$ be a differential
satisfying (5.7). For an even $a\in \A$ put
$$
\delta_a:=\delta +\p_a,\ \delta_a(b)=\d b+[a\b b].
 \eqno(5.10)
$$
Then we have $\d_a^2=0$ if the odd Maurer--Cartan equation is
satisfied:
$$
\d a+\frac{1}{2}[a\b a]=0.
 \eqno(5.11)
$$
Furthermore,
$$
[\d_a,\D]=0\quad \roman{if}\quad \D a =0.
 \eqno(5.12)
$$
Therefore,  from (5.11), (5.12) it
follows that $(\A,\D ,\d_a)$ is a differential GBV--algebra (dGBV).
In particular, 
$$
[\d_a ,\p_b]=\p_{\d_ab}.
 \eqno(5.13)
$$
We can in the same way shift $\D$. The essential difference
is that, as $\D$ itself, the shifted differential generally
will not be the derivation with respect to the
associative multiplication.

\medskip

{\bf 5.3. Homology of $(\A,\delta ).$} Since $\d$ is the derivation
with respect to both multiplications in $\A$
(associative one and the bracket), $\roman{Ker}\,\d$ is
the subalgebra with respect to both of them,
and $\roman{Im}\,\d$ is the ideal in this subalgebra with respect
to both structures. Therefore the homology group $H(\A,\d )$
inherits both multiplications, satisfying the identities
(5.6) and (5.7).

\smallskip

This reasoning holds for $H(\A,\d_a )$ as well, if $a$
satisfies the Maurer--Cartan equation (5.11).

\medskip

{\bf 5.4. Homology of $(\A,\Delta ).$} The same reasoning
furnishes  only the structure of
odd Lie algebra on $H(\A,\Delta )$, because $\D$ is not a derivation with respect
to the associative multiplication. However, if $\d$ and
$\D$ satisfy conditions (A) and (B) below, we will have
the natural isomorphism $H(\A,\Delta )=H(\A,\delta ).$

\medskip

The Lemma below is well known, see e.g. [GoM].

\proclaim{\quad 5.4.1. Lemma}
Let $\A$ be an additive group supplied with two endomorphisms
$\d$ and $\D$ satisfying $\d^2=\D^2=0$ and $\d\D=\alpha\D\d$
where $\alpha$ is an automorphism of $\A$ such that
$\alpha (\roman{Im}\,\D\d )= {Im}\,\D\d .$
Then clearly,
$\roman{Im}\,\d\D=\roman{Im}\,\D\d\subset
\roman{Im}\,\d\cap\roman{Ker}\,\D$ and similarly
with $\d$ and $\D$  permuted.

\smallskip

The following statements are equivalent:

\medskip

(i) The inclusions of the differential subgroups
 $i:\,(\roman{Ker}\,\D,\d )\subset (\A,\d)$
and $j:\,(\roman{Ker}\,\d,\D )\subset (\A,\D)$
are quasi--isomorphisms (that is, induce isomorphisms
of homology).

\medskip

(ii) We have actually equalities:
$$
\roman{Im}\,\d\D=\roman{Im}\,\D\d=
\roman{Im}\,\d\cap\roman{Ker}\,\D,
 \eqno(A)
$$
$$
\roman{Im}\,\d\D=\roman{Im}\,\D\d=
\roman{Im}\,\D\cap\roman{Ker}\,\d.
 \eqno(B)
$$

\medskip

Assume that these conditions are satisfied.
Then the both homology groups in (i) are naturally
isomorphic to 
$$
(\roman{Ker}\,\D\,\cap\,\roman{Ker}\,\d)/\roman{Im}\,\d\D.
$$
Moreover, the natural map $\roman{Ker}\,\D \to H(\A,\D )$
induces the surjection of the differential groups
$(\roman{Ker}\,\D,\d ) \to (H(\A,\D ),0)$ which
is a quasi--isomorphism, and similarly with
$\d$ and $\D$ interchanged. Hence the
both differential groups $(\A,\D )$ and $(\A,\d )$
are formal.
\endproclaim

\medskip

{\bf Proof.} We have:
$$
H(i)\ \roman{is\ injective}\ \Longleftrightarrow \ 
\roman{Ker}\,\D\,\cap\,\roman{Im}\,\d = \d\,(\roman{Ker}\,\D).
 \eqno(5.14)
$$
$$
H(i)\ \roman{is\ surjective}\ \Longleftrightarrow \ 
\roman{Ker}\,\d\,\subset\,\roman{Ker}\,\D +\roman{Im}\,\d 
\Rightarrow \D\,(\roman{Ker}\,\d )=\roman{Im}\,\d\D.
 \eqno(5.15)
$$
Here and below all kernel and images are taken in $\A$.
In the right hand side of (5.14), the inclusion
$\supset$ is evident, and the injectivity of $H(i)$
supplies the reverse inclusion.
The last arrow in (5.15) is obtained by applying
$\D$ to the previous inclusion: this gives
$\D\,(\roman{Ker}\,\d)\subset \D\,(\roman{Im}\,\d )=\roman{Im}\,\d\D$
whereas the reverse inclusion is obvious.

\smallskip 

Interchanging $\d$ and $\D$ we find
$$
H(j)\ \roman{is\ injective}\ \Longleftrightarrow \ 
\roman{Ker}\,\d\,\cap\,\roman{Im}\,\D = \D\,(\roman{Ker}\,\d).
 \eqno(5.16)
$$
$$
H(j)\ \roman{is\ surjective}\ \Longleftrightarrow \ 
\roman{Ker}\,\D\,\subset\,\roman{Ker}\,\d +\roman{Im}\,\D 
\Rightarrow \d\,(\roman{Ker}\,\D )=\roman{Im}\,\d\D.
 \eqno(5.17)
$$
Taken together, (5.14) and (5.17) prove (A), and (5.15) and (5.16)
prove (B), so that we have established the implication
(i) $\Rightarrow$ (ii). 

\medskip

Conversely, assume that (A) and (B) hold.

\smallskip

Then $H(i)$
induces surjection on the homology, because if
$\d a=0,$ we have $\D a\in \roman{Ker}\,\d\cap\roman{Im}\,\D$
so that by (B), $\D a =\D\d b,$ and then $a-\d b\in \roman{Ker}\,\D$
represents the same homology class as $a.$

\smallskip

Moreover, $H(i)$ induces injection on the homology, because
if $a\in \roman{Ker}\,\D,\ a=\d b$ for some $b\in \A,$
then $a\in \roman{Ker}\,\D \,\cap\,\roman{Im}\,\d$ so that
by (A), $a=\d c$ for some $c\in \roman{Im}\,\D \subset \roman{Ker}\,\D.$

\smallskip

By symmetry, the same holds for $H(j).$

\smallskip

The cycle subgroup for both differential groups ($\roman{Ker}\,\d 
,\D )$ and $(\roman{Ker}\,\D ,\d )$ is $(\roman{Ker}\,\d\, 
\cap\,\roman{Ker}\,\D )$, and if (i) and (ii) hold,
the boundaries can be identified with $\roman{Im}\,\d\D$,
cf. (5.15) and (5.17). It remains to deduce formality,
say, from (A) and (B).

\medskip

The natural map $\roman{Ker}\,\D \to H(\A,\D )$
is compatible with differentials, because if
$a\in \roman{Ker}\,\D ,$ then $\d a\in \roman{Im}\,\d \cap
\roman{Ker}\,\D$ so that by (A), $\d a=\D\d b$ for some $b$,
and hence the map is compatible with the zero differential
on $H(\A,\D )$.

\smallskip

This map is surjective on the homology. In fact, consider the
class of $a,\ \D a=0$ in  $H(\A,\D ).$ Then $\d a\in \roman{Im}\,\d\cap
\roman{Ker}\,\D$ so that in view of (A), $\d a=\d\D b$,
and the $\d$--cycle $a-\D b$ represents the same class as $a.$

\smallskip

Finally, the map is injective on the homology. In fact,
if $a\in \roman{Im}\,\D$ and $\d a=0$, then in view
of (B), $a\in\roman{Im}\,\d\D\subset \d\,(\roman{Ker}\,\D ).$

\smallskip

Thus we established the two--step quasi--isomorphism
of $(\A,\d )$ with $(H(\A,\D ),0)$ and by symmetry
of $(\A,\D )$ with $(H(\A,\d ),0).$ But the first two groups
are also naturally quasi--isomorphic. So they are formal.

\bigskip 

{\bf 5.4.2. Remarks.} In the context of
dGBV--algebras, we will apply this identification
 to $(\A,\D ,\d_a)$ with variable or formal generic 
$a$. Then we will be able
to interpret the ``constant'' space $H=H(\A,\D )$
as the flat structure on the family of algebras
$H(\A,\d_a )$ parametrized by the points of the generic
formal section of the
Maurer--Cartan manifold. The important technical problem
will be then deriving the conditions
(A) and (B) for the variable $a$.

\medskip

Notice that taken together, (A) and (B) are equivalent to
$$
\roman{Im}\,\d\D=\roman{Im}\,\D\d=
(\roman{Ker}\,\d\cap\roman{Ker}\,\D)\cap\,(\roman{Im}\,\d + \roman{Im}\,\D).
 \eqno(C)
$$
To deduce, say, (A) from (C), one omits the last term
in (C) and gets $\roman{Im}\,\d\D\supset \roman{Im}\,\d
\cap \roman{Ker}\,\D$ whereas the inverse inclusion is obvious.
Similarly, (C) follows from (A) and (B) together.

\smallskip

Assume that $\A$ is finite dimensional over a field and
$\d$ varies in a family, say $\{\d_a\}$. After a generalization,
dimension of $\roman{Im}\,\d\D$ can only jump,
and that of on $\roman{Ker}\,\d$ only drop. Hence if (B)
holds at a point, it holds in an open neighborhood of it.
In the case of the Dolbeault complex (cf
[BK]), only the cohomology
will be finite--dimensional. The validity of (C) for a particular
$\d =\d_0$ follows from the K\"ahler formalism. The argument
of the previous section (Theorem 4.2 (ii)) furnishes the same result
for the generic formal deformation.

\bigskip

{\bf 5.5. Integral.} Let $(\A,\D ,\d )$ be a 
dGBV--algebra. An even $k$--linear functional $\int :\A\to k$
is called {\it an integral} if the following two conditions are
satisfied:
$$
\forall a,b \in \A,\ \int (\d a)b=(-1)^{\widetilde{a}+1}\int a\d b,
 \eqno(5.18)
$$
$$
\forall a,b \in \A,\ \int (\D a)b=(-1)^{\widetilde{a}}\int a\D b.
 \eqno(5.19)
$$
Notice that (5.18) is equivalent to $\forall a\in \A, \int \d a =0$
because $\d$ is a $k$--derivation. Applying (5.19) to $b=1$,
we see that $\forall a\in \A, \int \D a =0$ as well.

\medskip

\proclaim{\quad 5.5.1. Proposition}  Let $\int$ be an integral
for $(\A,\D ,\d ).$ 

\smallskip

(i) If $a$ or $b$ belongs to $\roman{Ker}\,\D,$ we have
$$
 \int \p_ab=\int [a\b b]=0.
 \eqno(5.20)
$$
Hence if $a$ satisfies (5.11) and (5.12), $\int$ is an integral
for $(\A,\D ,\d_a )$ as well.

\smallskip

(ii) $\int$ induces a linear functional on $H(\A,\D )$ and
$H(\A,\d_a)$ for all $a$ as above. These functionals
are compatible with the identifications following from
the condition (C).
\endproclaim

\smallskip

{\bf Proof.} If, say, $\D a=0,$ we have
$$
\int \p_ab= \int ((-1)^{\widetilde{a}}\D (ab)-(-1)^{\widetilde{a}}(\D a)\,b-
a\,\D b)= -\int a\,\D b=-(-1)^{\widetilde{a}}\int\D a\,b=0.
$$
The rest is straightforward.

\medskip

{\bf 5.6. Metric.} If $\int$ is an integral on $(\A,\D ,\d )$,
we can define the scalar products on $H(\A,\d_a)$ induced
by the symmetric scalar product $(a,b)\mapsto \int ab$ on $\A.$
For the construction of Frobenius manifolds,
it is necessary to ensure that these scalar products
are non--degenerate.

\smallskip

Integral and metric are compatible with base extensions.

\medskip

{\bf 5.7. Additional grading.} Assume now that $\A$
as commutative $k$--superalgebra is graded
by an additive subgroup of $k$. Thus $\A=\oplus_n\A^n,$
$k\in \A^0$, $\A^m\A^n\subset \A^{m+n},$ and each
$\A^i$ is graded by parity. We write $|a|=i$ if $a\in \A^i.$
Various induced gradings and degrees of homogeneous operations
are denoted in the same way. (In the main example of [BK],
$\A$ is $\bold{Z}$--graded, and each $\A^i$ is either even, or odd,
but this plays no role in general).

\smallskip

All base extensions then must be furnished by the
similar grading or its topological completion.

\smallskip

We will assume also that $|\D |=-1.$ It follows 
that $|[a\b b]|=|a|+|b|-1$ which we interpret as 
$|\b |=-1.$ Moreover, we postulate that $|\d |=1.$
This means that the shifted differential $\d_{\gamma}$
can be homogeneous only for $|\gamma |=2,$
and similarly for extended base.

\smallskip

Homology space $H$ in all its incarnations
(cf. Lemma 5.4.1) inherits the grading from $\A .$
The dual space $H^t$ is graded in such a way that
the pairing $H^t\otimes H\to k$ has degree zero.
This induces the additional grading (or more precisely,
the notion of homogeneity) on $K=k[[H^t]]$ (which might
be the product rather than the sum of its homogeneous
components).

\smallskip

Integral is supposed to have a definite degree,
not necessarily zero (and usually non--zero).
 
\medskip

{\bf 5.8. Tensor product of GBV--algebras.} Let $(\A_i,\D_i),\,i=1,2,$
be two GBV--algebras over $k$. Put $\A:=\A_1\otimes \A_2,
\D :=\D_1\otimes 1+1\otimes \D_2:\,\A\to \A.$ 

\medskip
\proclaim{\quad 5.8.1. Proposition} $(\A,\D )$ is a GBV--algebra.
We have for $a_i,b_i\in \A_i$
$$
\p_{a_1\otimes a_2}=\p_{a_1}\otimes (-1)^{\widetilde{a}_2}l_{a_2}+
l_{a_1}\otimes \p_{a_2},
 \eqno(5.21)
$$
or equivalently
$$
[a_1\otimes a_2\bullet b_1\otimes b_2]=
(-1)^{\widetilde{a}_2(\widetilde{b}_1+1)}
[a_1\bullet b_1]\otimes a_2b_2+
(-1)^{\widetilde{b}_1(\widetilde{a}_2+1)}
a_1 b_1\otimes [a_2\bullet b_2].
 \eqno(5.22)
$$
\endproclaim

\smallskip

{\bf Proof.} (5.21) is 
established by a straightforward calculation which we omit.
From (5.21) it follows that $\p_{a_1\otimes a_2}$ are derivations.
Hence $\p_a$ are derivations for all $a\in \A$ so that $(\A,\D )$
is a GBV--algebra. (5.22) is a rewriting of (5.21).

\medskip

Clearly, tensor product is commutative and associative with respect to the
standard isomorphisms.

\medskip

If $\d_i:\,\A_i\to \A_i$ are odd derivations of $(\A_i,\D_i)$
satisfying (5.7), then $\d :=\d_1\otimes 1+1\otimes \d_2$
is an odd derivation of $\A_1\otimes \A_2$ satisfying (5.7).

\smallskip

If $\A_i$ are furnished with additional gradings having
the properties postulated above, then the total
grading on $\A_1\otimes \A_2$ satisfies the same conditions.

\bigskip

{\bf 5.8.2. Decomposable solutions to the Maurer--Cartan equation.}
In the notation of the previous subsection,
let $(\A,\D ,\d )$ be the tensor product of $(\A_i,\D_i ,\d_i ),
i=1,2.$ Assume that $a_i\in \A_i$ satisfy the Maurer--Cartan
equation (5.11). Then from (5.22) it follows that
$a:=a_1\otimes 1+1\otimes a_2$
satisfies (5.11) as well.  Moreover, if $\D_ia_i=0$, then
$\D a=0,$ so that $(\A,\D ,\d_a)$ is the differential
GBV--algebra. Such structures will be called
decomposable ones.

\newpage

\centerline{\bf \S 6. From dGBV--algebras to Frobenius manifolds}

\bigskip

{\bf 6.1. Normalized formal solution to the master equation.}
In this section, we fix a dGBV $k$--algebra $(\A ,\D ,\d )$
and the derived odd bracket $[\b ]$ on it. We will assume
that this algebra satisfies a series of assumptions
which will be introduced and numbered consecutively.

\medskip

{\it ASSUMPTION 1. $(\A ,\D ,\d )$ satisfies conditions
of the Lemma 5.4.1. Moreover, the homology group
$H=H(\A ,\d )$ (and any group naturally isomorphic
to it) is a free $k$--module of finite rank.}

\medskip

Choosing an indexed basis $[c_i], c_i\in \A$ of $H$ 
and the dual basis $(x_i)$ of $H^t$ we will always
assume that $c_0=1.$ 
As in 4.1, but now conserving parity, we put $K:=k[[H^t]]=
k[[x_i]].$
We will denote by $X_i=\p /\p x_i$ the respective partial
derivatives acting on $K$ and on $K\otimes \A , K\otimes H$
etc via the first factor.

\smallskip

\proclaim{\quad 6.1.1. Proposition} If $(\A,\D ,\d )$
satisfies Assumption 1 above, then there exists
a generic even formal solution  $\Gamma =\sum_i\Gamma_i\in K\otimes 
\roman{Ker}\,\D$ to the master
equation 
$$
\d \Gamma +\frac{1}{2}[\Gamma \b\Gamma]=0
$$
with the following properties:

\smallskip

(i) $\Gamma_0=0, \Gamma_1=\sum x_ic_i, \Gamma_n\in K\otimes
\roman{Im}\,\D$ for all $n\ge 2.$ Here $c_i\in 
\roman{Ker}\,\D\,\cap\,\roman{Ker}\,\d$, 
and $\Gamma_n$ is the homogeneous component 
of $\Gamma$ of degree $n$ in $(x_i).$

\smallskip

(ii) Moreover, this $\Gamma$ can be chosen
in such a way that $X_0\Gamma = 1.$

\smallskip

Such a solution will be called normalized.

\endproclaim
 
\smallskip

{\bf Proof.} The first statement follows from
the Theorem 4.2 (i) applied to the odd differential Lie superalgebra
$(\roman{Ker}\,\D , [\b ],\d ).$ 

\smallskip

We must only check that the conditions of the applicability of
this theorem are satisfied. To facilitate the bookkeeping
for the reader, we register the correspondences between the old
and the new notation:
$\roman{Ker}\,d$ becomes $\roman{Ker}\,\D\,\cap
\roman{Ker}\,\d,$ $\roman{Ker}\,\phi$ turns into
$\roman{Im}\,\D,$ $\roman{Im}\, d$ corresponds to
$\roman{Im}\,\d\D.$ All of this forms a part of Lemma
5.4.1. From the second formula in 5.1 it follows
that $\roman{Ker}\,\D$ is closed with respect to $[\b ]:$
$[a\b b]=(-1)^{\widetilde{a}}\D (ab).$ This formula
shows as well that $[\b ]$ induces zero operation on
$H(\roman{Ker}\,\D,\d):$ if $a,b\in \roman{Ker}\,\D\,\cap
\roman{Ker}\,\d,$ then $[a\b b]\in \roman{Im}\,\D\,\cap
\roman{Ker}\,\d = \roman{Im}\,\d\D.$

\smallskip

It remains to check the assertion (ii). Clearly,
our choice $c_0=1$ assures that $X_0\Gamma_1=1.$
Assume by induction that $\Gamma_2 ,\ldots ,\Gamma_n$
do not depend on $x_0$. Clearly,
$[c_0\b a]$=0 for any $a,$ so that
in the (odd version of the) equation (4.5)
the right hand side is independent of $x_0$ as well.
Since the argument showing the existence
of $\Gamma_{n+1}$ in the proof of 4.2 can be
applied to each coefficient of the monomials in $x_i$
separately, we may find $\Gamma_{n+1}$ independent of
$x_0.$ The final normalization argument
can be also applied coefficientwise.

\medskip

{\bf 6.2. The (pre)--Frobenius manifold associated to $(\A ,\D ,\d ).$}
Consider the formal manifold $M$, the formal spectrum of
$K$ over $k.$ The flat coordinates will be by definition
$(x_i)$ so that the space of flat vector fields can be
canonically identified with $H.$ We fix a normalized $\Gamma$ as above. 

\medskip

\proclaim{\quad 6.2.1. Lemma} The bi--differential group
$(\A_K ,\D_K ,\d_{\Gamma})$ satisfies the conditions
and conclusions of Lemma 5.4.1.
\endproclaim

\smallskip

{\bf Proof.} Clearly, $\d_{\Gamma}^2=\D_K^2=[\D_K ,\d_{\Gamma}]=0$
(the latter follows from (5.3)). From the Assumption 1
and the proof of 6.6.1 above we see that we can
apply the Theorem 4.2 (ii) to $(K\otimes \roman{Ker}\,\D,\d_{\Gamma})$
instead of $(g_K,d_{K,\Gamma}).$ The inclusion (4.6)
reads in this context:
$$
\roman{Im}\,\D_K\,\cap\,\roman{Ker}\,\d_{\Gamma}\subset
\roman{Im}\,\d_{\Gamma}\D_K
$$
which implies the condition (B).

\smallskip

To check (A), consider the inclusion map
$$
\roman{Im}\,\d_{\Gamma}\D_K\to \roman{Im}\,\d_{\Gamma}\,
\cap\,\roman{Ker}\,\D_K .
$$
It becomes an isomorphism after reduction modulo $(x_i)$
in view of the Assumption 1. Hence it is an isomorphism.

\medskip

{\bf 6.2.2. $\circ$--multiplication on tangent fields.}
We will define now
the $K$--linear $\circ$--multiplication on the $K$--module
of all vector fields $\Cal{T}_M=K\otimes H=H_K.$

\smallskip

To this end we first apply Theorem 4.2 (ii)   to
the odd differential Lie algebra $(\roman{Ker}\,\D ,\d_{\Gamma}).$
It shows that the homology of this
algebra is naturally identified with $H_K.$

\smallskip

In view of Lemmas 5.4.1 and 6.2.1 we know that
the injection $(\roman{Ker}\,\D_K ,\d_{\Gamma})\to
(\A_K,\d_{\Gamma})$ induces isomorphism of homology
$H_K= \roman{Ker}\,\d_{\Gamma}/\roman{Im}\,\d_{\Gamma}.$
But $\roman{Ker}\,\d_{\Gamma}$ is a commutative
$K$--subalgebra of $A_K$ and $\roman{Im}\,\d_{\Gamma}$
is an ideal in it. Hence $H_K$ inherits the
multiplication which we denote $\circ$. We record
the following ``explicit'' formula for it.
Interpreting any $X\in H_K$ as the derivation
$\overline{X}$ of $K\otimes \A$ acting through
the first factor (cf. Corollary 4.2.1), we have:
$$
\overline{X\circ Y}\Gamma \equiv
\overline{X}\Gamma\cdot\overline{Y}\Gamma\,
\roman{mod}\,\roman{Im}\,\d_{\Gamma}
\eqno(6.1)
$$
(dot here means the associative multiplication in $\A_K$). 
This follows directly from the Corollary 4.2.1
applied to our situation. Notice that whereas
$\overline{X}\Gamma$ and $\overline{Y}\Gamma$
lie in $\roman{Ker}\,\D_K\,\cap\,\roman{Ker}\,\d_{\Gamma},$
their product generally lies only in the larger group
$\roman{Ker}\,\d_{\Gamma}.$

\smallskip

Directly from the initial definition one sees that
$e:=X_0$ is the flat identity for $\circ .$

\smallskip

In order to complete the description of the pre--Frobenius
structure, it remains to choose a flat metric on $M$.

\medskip

{\it ASSUMPTION 2. There exists an integral $\int$
for $(\A ,\D ,\d )$ such that the bilinear
form on $H=H(A,\d )$ induced by
$(X,Y)\mapsto \int \overline{X}\Gamma \cdot\overline{Y}\Gamma$
is non--degenerate.}

\medskip

Denoting this form $g$ we clearly have the invariance
property defining the symmetric multiplication tensor $A$:
$$
g(X,Y\circ Z)=g(X\circ Y,Z):= A(X,Y,Z).
\eqno(6.2)
$$
We will check now that this structure is actually
Frobenius.
Since the $\circ$--multiplication is associative,
we have only to establish its potentiality.

\smallskip

To this end we will check Dubrovin's criterium:
the structure connection $\nabla_{\lambda}$
on $\Cal{T}_M$ is flat (cf. [M], Ch. I, Theorem 1.5).

\smallskip

To be more precise, let $\nabla_0$ be the flat
connection on $\Cal{T}_M$ whose horizontal
sections are $H$. Clearly, $\nabla_{0,Y}(Z)=
\overline{Y}(Z)$ where this time $\overline{Y}$
means $Y$ acting on $K\otimes H$ via $K.$
By definition, 
$$
\nabla_{\lambda ,Y}(Z)=\overline{Y}(Z)+\lambda Y\circ Z
\eqno(6.3)
$$
where $\lambda$ is an even parameter.

\smallskip

We have the canonical surjection $\roman{Ker}\,\D_K
\to H_K$ and the two lifts of $X$ both denoted by
$\overline{X}$ are compatible with this surjection,
and also with embedding $\roman{Ker}\,\D_K\subset \A_K.$
Therefore the section of $\roman{Ker}\,\D_K\to H_K$
denoted $\psi$ in the Corollary 4.2.1 sends
$\overline{Y}(Z)$ to $\overline{Y}(\overline{Z}\Gamma ),$
and $\nabla_{\lambda ,Y}(Z)$ lifts to $\overline{Y}(\overline{Z}\Gamma )
+\lambda\overline{Y}\Gamma\cdot\overline{Z}\Gamma$
in view of (6.1) and (6.3).

\smallskip

Our preparations being now completed, we can prove

\medskip

\proclaim{\quad 6.2.3. Theorem} The connection $\nabla_{\lambda}$
is flat. Hence the pre--Frobenius structure defined above is
potential.
\endproclaim

\smallskip

{\bf Proof.} Applying (6.3) twice, we find
$$
\nabla_{\lambda ,X}\nabla_{\lambda ,Y}(Z)=
\overline{X}(\overline{Y}(Z))+\lambda\overline{X} (Y\circ Z)+
\lambda X\circ\overline{Y} (Z)+\lambda^2 X\circ Y\circ Z.
\eqno(6.4)
$$
We may and will consider only the case when
$X,Y$ supercommute (e.g. $X,Y\in H$).
In order to establish flatness, it suffices to check that
$$
\overline{X} (Y\circ Z)+X\circ\overline{Y} (Z)=
(-1)^{\widetilde{X}\widetilde{Y}}
(\overline{Y} (X\circ Z)+Y\circ\overline{X} (Z)).
\eqno(6.5)
$$
We will see that already the
$\psi$--lifts of both sides of (6.5) coincide up to $\roman{Im}\,\d_{\Gamma}$.
In fact, $\overline{X} (Y\circ Z)$ lifts to 
$\overline{X}(\overline{Y}\Gamma\cdot\overline{Z}\Gamma ),$
$X\circ\overline{Y} (Z)$ lifts to 
$\overline{X}\Gamma\cdot\overline{Y}(\overline{Z}\Gamma )$
so that (6.5) becomes
$$
\overline{X}(\overline{Y}\Gamma )\cdot\overline{Z}\Gamma
+(-1)^{\widetilde{X}\widetilde{Y}}
\overline{Y}\Gamma\cdot\overline{X}(\overline{Z}\Gamma )
+\overline{X}\Gamma\cdot\overline{Y}(\overline{Z}\Gamma )=
$$
$$
(-1)^{\widetilde{X}\widetilde{Y}}
\overline{Y}(\overline{X}\Gamma )\cdot\overline{Z}\Gamma
+\overline{X}\Gamma\cdot\overline{Y}(\overline{Z}\Gamma )+
(-1)^{\widetilde{X}\widetilde{Y}}
\overline{Y}\Gamma\cdot\overline{X}(\overline{Z}\Gamma ).
$$
This finishes the proof.

\medskip

{\bf 6.3. Euler field.} Assume now that $\A$
is endowed with a grading satisfying the
conditions of 5.7. All the previous discussion
makes sense, and the results hold true, if we add appropriate
grading conditions at certain places, the
most important of which is $|\Gamma |=2$
implying $|\circ |=2$ in view of (6.1).

\smallskip

Denote by $E$ the derivation of $K$ defined by the
following Euler condition:
$$
\forall f\in K,\ Ef=\frac{1}{2}|f|f,
\eqno(6.6)
$$
where $|\,\,|$ is the grading induced on $H^t$
from $\A$ via $H.$ In cooordinates as in 6.1 we have
$$
E=\frac{1}{2}\sum_i|X_i|\,x_iX_i.
\eqno(6.7)
$$

{\it ASSUMPTION 3. Assume that the integral
is homogeneous and denote its degree by $-2D-4.$}

\smallskip

For the general discussion of spectrum
cf. [M], Chapter I, \S 2, cf. also 1.3 above for
spectrum of quantum cohomology. 

\medskip

\proclaim{\quad 6.3.1. Proposition} $E$ is an Euler
field on the formal Frobenius manifold
described in 6.2. Its spectrum is
$(D;\, d_i$ with multiplicity $\roman{dim}\, H^{-2d_i}$),
and $d_0=1.$
\endproclaim

\smallskip

{\bf Proof.} Comparing
(6.7) with the notation of [M], we see
that the spectrum of $-\roman{ad}\, E$ on $H=\Cal{T}_M^f$
is $d_i$ with multiplicity $\roman{dim}\, H^{-2d_i}$
where 
$$
d_i=\frac{1}{2}\,|x_i|=-\frac{1}{2}\,|X_i|.
\eqno(6.8)
$$

\smallskip

Since $X_0\circ X=X$ and $|\circ |=2$, we have $d_0=1.$

We must now check the formula
$$
E(g(X,Y))-g([E,X],Y)-g(X,[E,Y])=Dg(X,Y).
\eqno(6.9)
$$
It suffices to do this for the case when $X,Y$ are flat vector fields
having definite degrees. Then $[E,X]= -\dfrac{|X|}{2}\,X.$
 Since $g(X,Y)\in k,$
(6.9) becomes
$$
(|X|+|Y|-2D)g(X,Y)=0.
$$ 
But  $g(X,Y)=\int \overline{X}\Gamma
\cdot \overline{Y}\Gamma$ vanishes unless $-2D-4+|X|+2+|Y|+2=0$
which proves (6.9). 

\smallskip

Furthermore, from (6.1) we infer that $|X\circ Y|=|X|+|Y|+2.$
Hence if $X_i\circ X_j=\sum_kA_{ij}{}^kX_k,$ we have
$$
E\,A_{ij}{}^k=\frac{1}{2}\,|A_{ij}{}^k|\cdot A_{ij}{}^k=
\frac{1}{2}\,(|X_i|+|X_j|-|X_k|+2)\,A_{ij}{}^k.
$$
Comparing this with the formula (2.18) of [M], Chapter I
and taking into account (6.8), 
we see that $E$ satisfies the Definition 2.2.1 of
[M], loc. cit. This finishes the proof.

\medskip

Notice that the Euler field (6.7) contains no flat
summand: $X_i$ with $d_i=0$ do not contribute.
Hence if this construction furnishes a Frobenius
manifold which is quantum cohomology of some $V$,
then $c_1(V)$ must vanish (modulo torsion).
 
\medskip

{\bf 6.4. Explicit potential.}
The direct way to establish potentiality 
is to find an even series $\Phi \in K$
such that  for all $X,Y,Z\in H$ we have 
$A(X,Y,Z)=XYZ\Phi $ (from now on, we write
$X$ instead of $\overline{X}$ in order to denote
derivations on various $K$--modules acting through $K$).
Moreover, it suffices to check this for $X=Y=Z.$
We will give here the beautiful formula 
of Chern--Simons type for $\Phi$ discovered in [BK].

\smallskip

Extend the integral to the $K$--linear map $\int :\,\A_K \to K.$
For a fixed normalized $\Gamma$ put $\Gamma =\Gamma_1 +\D B$
where $B_0=B_1=0$ and $\Delta$ means $\D_K.$ 

\smallskip

\proclaim{\quad 6.4.1. Theorem} The formal function  
$$
\Phi := \int \left(\frac{1}{6}\, \Gamma^3-\frac{1}{2}\,\d B\,\D B\right) .
\eqno(6.10)
$$
is a potential for the Frobenius manifold defined above.
\endproclaim

\smallskip

{\bf Proof.} We have to prove that for
any $X\in H$
$$
A(X,X,X)=\int (X\Gamma )^3=X^3\Phi .
\eqno(6.11)
$$

\smallskip

We supply below the detailed calculation
consisting of the series of elementary steps,
each being an application of one
of the identities (5.18), (5.19), Leibnitz
rule for (super)derivations and the fact
that $\d ,\D , X$ pairwise supercommute.
Moreover, we use the master equation in the form
$\D\Gamma^2=-2\,\d\Gamma$ following from
$\D\Gamma =0.$ Finally, $\d\Gamma_1=X^n\Gamma_1=0$
for $n\ge 2$ so that $\d\Gamma =\d\D B, X^n\Gamma =X^n\D B.$

\smallskip

We start with treating the first summand of the right hand side of (6.10).
The derivation $X$ is interchangeable
with integration, so we have by the Leibnitz rule
$$
X^3\left( \frac{1}{6}\int \Gamma^3\right)=
\int  (X\Gamma)^3 +\int \left((2+(-1)^{\widetilde{X}})\Gamma\cdot X\Gamma\cdot
X^2\Gamma +\frac{1}{2}\,\Gamma^2\cdot X^3\Gamma \right) .
\eqno(6.12)
$$
The second summand of (6.10) is added in order to
cancel the extra terms in (6.12). 

\smallskip

First, we rewrite it:
$$
\frac{1}{2}\int \d B\,\D B=\frac{1}{2}\int  B\,\d\D B=
\frac{1}{2}\int B\,\d \Gamma = -\frac{1}{4}\int B\,\D(\Gamma^2)=
\frac{1}{4}\int \D B\cdot\Gamma^2.
\eqno(6.13)
$$
(We could have chosen the last expression in (6.13) from the start).

\smallskip
Now, again by Leibnitz rule,
$$
\frac{1}{4}\,X^3\int \D B\cdot\Gamma^2=
\frac{1}{4}\int \left( X^3(\D B)\,\Gamma^2
+(2+(-1)^{\widetilde{X}})X^2(\D B)\cdot X(\Gamma^2)\right .
$$
$$
\left.  +(2+(-1)^{\widetilde{X}})X(\D B)\cdot X^2(\Gamma^2)+
\D B\cdot X^3(\Gamma^2)\right) .
\eqno(6.14)
$$
The first two summands in (6.14) can be 
directly rewritten in the same form as
in (6.12):
$$
 X^3(\D B)\cdot \Gamma^2
+(2+(-1)^{\widetilde{X}})X^2(\D B)\cdot X(\Gamma^2)=
$$
$$
 X^3\Gamma\cdot \Gamma^2
+(2+(-1)^{\widetilde{X}})X^2\Gamma\cdot X(\Gamma^2).
\eqno(6.15)
$$
The third summand takes somewhat more work:
$$
\int X(\D B)\cdot X^2(\Gamma^2)= 
-\int X B\cdot X^2(\D\Gamma^2)= 
2\int X B\cdot X^2(\d\Gamma )=
$$
$$
2\int X B\cdot \d (X^2\Gamma )=
2\int X B\cdot \d \D(X^2B )=
-2\int \d \D X B\cdot X^2B=
$$
$$
-2\int X\d \Gamma \cdot X^2B=
\int X\D (\Gamma^2) \cdot X^2B=
\int X(\Gamma^2)\cdot X^2\D B=
$$
$$
\int X(\Gamma^2)\cdot X^2\Gamma .
\eqno(6.16)
$$
\medskip

Finally, the fourth summand is calculated similarly,
but in two steps. We start with an expression
of the second order in $X$:
$$
\int \D B\cdot X^2(\Gamma^2)= 
-\int  B\cdot X^2(\D\Gamma^2)= 
2\int  B\cdot X^2(\d\Gamma )=
$$
$$
2\int  B\cdot \d (X^2\Gamma )=
2\int  B\cdot \d \D(X^2B )=
-2\int \d \D B\cdot X^2B=
$$
$$
-2\int \d \Gamma \cdot X^2B=
\int \D (\Gamma^2) \cdot X^2B=
\int \Gamma^2\cdot X^2\D B=
$$
$$
\int \Gamma^2\cdot X^2\Gamma .
\eqno(6.17)
$$ 
Apply now $X$ to the first and the last expressions
of (6.17). We get
$$
\int X(\D B)\cdot X^2(\Gamma^2)+
\int \D B\cdot X^3(\Gamma^2)=
\int X(\Gamma^2)\cdot X^2\Gamma+
\int \Gamma^2\cdot X^3\Gamma .
\eqno(6.18)
$$
Comparing this with (6.16), one gets
$$
\int \D B\cdot X^3(\Gamma^2)=\int \Gamma^2\cdot X^3\Gamma .
\eqno(6.19)
$$
Putting all of this together,
one obtains finally (6.11).

\newpage

\centerline{\bf References}

\bigskip

[AGV] V.~Arnold, S.~Gusein--Zade, A.~Varchenko. {Singularities
of differentiable maps, vols. I, II.} Birkh\"auser, Boston,
!985 and 1988.

\smallskip

[BK] S.~Barannikov, M.~Kontsevich. {\it Frobenius manifolds and
formality of Lie algebras of polyvector fields.}
Preprint alg--geom/97010072

\smallskip

[B] K.~Behrend. {\it Gromov--Witten invariants in algebraic geometry.}
Inv. Math., 127 (1997), 601--617.

\smallskip

[BF] K.~Behrend, B.~Fantechi. {\it The intrinsic normal cone.}
Inv. Math., 128 (1997), 45--88.

\smallskip

[BM] K.~Behrend, Yu.~Manin. {\it Stacks of stable maps and
Gromov--Witten invariants.} Duke Math. J., 85:1 (1996), 1--60.

\smallskip

[C] S.~Cecotti. {\it $N=2$ Landau--Ginzburg vs. Calabi--Yau
$\sigma$--models: non--perturbative aspects.} Int. J. of Mod. Phys. A,
6:10 (1991), 1749--1813.

\smallskip

[D] B.~Dubrovin. {\it Geometry of 2D topological field theories.}
In: Springer LNM, 1620 (1996), 120--348

\smallskip

[DZh] B.~Dubrovin, Youjin Zhang. {\it Extended affine Weyl groups
and Frobenius manifolds.} Preprint SISSA 67/96/FM.

\smallskip

[Ge1] D.~Gepner. {\it On the spectrum of 2D conformal field theory.}
Nucl. Phys., B287 (1987), 111--126.

\smallskip

[Ge2] D.~Gepner. {\it Fusion rings and geometry.} Comm. Math. Phys.,
141 (1991), 381--411.

\smallskip

[Ger] M.~Gerstenhaber. {\it The cohomology structure of an
associative ring.} Ann. of Math., 79:1, 267--

\smallskip

[Giv] A.~Givental. {\it Equivariant Gromov--Witten invariants.}
Int. Math. Res. Notes, 13 (1996), 613--663.

\smallskip

[GoM] W.~Goldman, J.~Millson. {\it The deformation theory
of representations of fundamental groups of compact
K\"ahler manifolds.} Publ. Math. IHES, 86 (1988), 43--96.

\smallskip

[Gr] B.~Green. {\it Constructing mirror manifolds.}
In: Mirror Symmetry II,
ed.~by B.~Greene and S.~T.~Yau, AMS--International Press, 1996, 29--69.

\smallskip

[G--ZV] S.~Gusein--Zade, A.~Varchenko. {\it Verlinde algebras and
the intersection form on vanishing cycles.} Selecta Math.,
New. Ser., 3 (1997), 79--97.

\smallskip 

[K] R.Kaufmann, {\it The geometry of moduli spaces of pointed
curves, the tensor product in the theory of Frobenius manifolds,
and the explicit K\"unneth formula in quantum cohomology.}
Ph. D. thesis, MPI f\"ur Mathematik, Bonn, 1997.

\smallskip

[Ko] M.~Kontsevich. {\it Enumeration of rational curves via torus actions.}
In: The Moduli Space of Curves, ed. by
R\. Dijkgraaf, C\. Faber, G\. van der Geer, Progress in Math\.
vol\. 129, Birkh\"auser, 1995, 335--368.

\smallskip

[KM1] M.~Kontsevich, Yu.~Manin. {\it Gromov-Witten classes, quantum
cohomology, and enumerative geometry.} Comm. Math. Phys.,
164:3 (1994), 525--562.

\smallskip

[KM2] M.~Kontsevich, Yu.~Manin. {\it Relations between the correlators
of the topological sigma--model coupled to gravity.}
Preprint alg--geom/970824

\smallskip

[KMK] M.~Kontsevich, Yu.~Manin. {\it Quantum cohomology of a product (with  Appendix by R. Kaufmann)}. Inv. Math., 124, f. 1--3 (1996), 313--339.

\smallskip

[LZ] B.~H.~Lian, G.~Zuckerman. {\it New perspectives
on the BRST--algebraic structure of string theory.}
Comm. Math. Phys., 154 (1993), 613--646, hep--th/9211072

\smallskip

[M] Yu.~Manin. {\it Frobenius manifolds, quantum cohomology,
and moduli spaces (Chapters I,\,II,\,III).} Preprint MPI  96--113,
1996.

\smallskip

[MM] Yu.~Manin, S.~Merkulov. {\it Semisimple Frobenius
(super)manifolds and quantum cohomology of $\bold{P}^r$.}
Topological Methods in Nonlinear Analysis, 9:1 (1997), 107--161,
alg--geom/9702014.

\smallskip

[O] T.~Oda. {\it K.~Saito's period map for holomorphic functions with
isolated critical points.} In: Adv. Studies in Pure Math., 10 (1997),
Algebraic Geometry, Sendai, 1985, pp. 591--648.

\smallskip

[S1] K.~Saito. {\it Primitive forms for a universal unfolding
of a function with an isolated critical point.} Journ. Fac. Sci.
Univ. Tokyo, Sec. IA, 28:3 (1982), 775--792.

\smallskip

[S2] K.~Saito. {\it Period mapping associated to a primitive form.}
Publ. RIMS, Kyoto Univ., 19 (1983), 1231--1264.

\smallskip

[St] J.~Stasheff. {\it Deformation theory and the Batalin--Vilkovisky
master equation.} In: Deformation Theory
and Symplectic Geometry, Eds. D.~Sternheimer et al.,
Kluwer, 1997, 271--284.

\smallskip

[Z] J.-B.~Zuber. {\it Graphs and reflection groups.}
Comm. Math. Phys., 179 (1996), 265--294.

\enddocument